\documentclass[a4paper]{article}
\usepackage{multirow}

\usepackage{graphicx} 
\usepackage{mathtools,bm}
\usepackage{bbold}   
\usepackage{amsmath}
\usepackage{amsfonts}
\usepackage{courier}   
\usepackage{color}
\usepackage{tabu}
\usepackage{soul}
\usepackage{algorithm}
\usepackage[noend]{algpseudocode}
\usepackage{float}
\usepackage[top=1in, bottom=1in, left=1in, right=1in]{geometry}
\usepackage[titletoc,title]{appendix}
\usepackage{amsthm}
\usepackage{amssymb}
\usepackage{arydshln}
\usepackage[colorlinks,
linkcolor=blue,
anchorcolor=blue,
citecolor=blue]{hyperref}
\usepackage{verbatim}

\numberwithin{equation}{section}

\newcommand{\bs}{\boldsymbol}

\definecolor{lightblue}{rgb}{0.0,0.5,0.8}

\definecolor{darkgreen}{rgb}{0.0, 0.4, 0.2}
\newtheorem{theorem}{Theorem}[section]
\newtheorem{lemma}{Lemma}[section]

\makeatother

\title{A duality-based approach for solving linear parabolic control constrained optimal control problems}
\author{Hailing Wang}
\begin{document}
	\maketitle
	\section*{Abstract}
     	This paper is concerned with the optimal control problem governed by a linear parabolic equation and subjected to box constraints on control variables. This type of problem has important applications in heating and cooling systems. By applying the scheme of Fenchel duality, we derive the dual problem explicitly where the control constraints in primal problem are embedded in the dual problem's objective functional. The existence and uniqueness of the solution to the dual problem are proved and the first-order optimality conditions are also derived. In addition, we discuss the saddle point property between solution of the primal problem and the dual problem. The solution of primal problem can be readily obtained by the solution of the dual problem. To solve the dual problem numerically, we design two implementable method: conjugate gradient method and semismooth Newton method. Three example problems are solved, numerical results show that the proposed method is efficient and accurate.

	\section{Introduction}
	We consider the following constrained optimal control problem
	\begin{equation}\label{eq3}
	\min\limits_{ u\in\mathcal{C}} \frac{1}{2} \iint_{Q}\left|y-y_{d}\right|^{2} d x d t+\frac{\gamma}{2} \iint_{\mathcal{O}}|u|^{2} dx dt
	\end{equation}
	subject to
	\begin{equation}\label{eq2}
	\left\{
	\begin{aligned}
	&\frac{\partial y}{\partial t}-\nu\Delta y+a_{0}(x,t) y=u\chi_{\mathcal{O}}, &\text { in } \Omega \times(0, T), \\
	&y=0, & \text { on } \Gamma \times(0, T), \\
	&y(0)=y_0, & \text { in } \Omega.
	\end{aligned}
	\right.
	\end{equation}
Here, $Q=\Omega \times(0, T), \mathcal{O}=O\times(0,T)$ and $0<T<+\infty$, $\Omega \subset \mathbb{R}^{n}$ is the domain of space variable $x$ and $O$ is an open subset of $\Omega$, $\Gamma=\partial\Omega$ is the piecewise continuous boundary of $\Omega$. $u$ is the control variable, $\mathcal{C}$ is the admissible control set specified by
	\begin{equation}\label{eq4}
	\mathcal{C}=\left\{u \mid u \in L^{2}(\mathcal{O}), a \leq u(x , t) \leq b, \text { a.e. in } \mathcal{O}\right\}
	\end{equation}
$a$ and $b$ are given constants satisfying $a\leq 0\leq b$. $y$ is the solution of state equation (\ref{eq2}). $y_{d}$ is the target function given in $L^{2}(Q)$ and  
$\gamma>0$ is a regularization parameter. $\nu$ is a positive constant and $a_0(x,t)\geq 0$ is a function given in functional space $L^\infty(Q)$, $\chi_{\mathcal{O}}$ denotes the characteristic  function of set $\mathcal{O}$ and $y_0\in L^2(\Omega)$ is the initial state. The existence and uniqueness of the optimal control to this problem can be found in \cite{MR0271512}.
	
Such optimal control problem has wide applications in heating and cooling systems \cite{MR1860629,appl,article,doi:10.1142/S0218202509004054}. Since it is impossible to obtain an analytic solution for this problem. Numerical methods are indispensable.

Notice that the problem (\ref{eq3})-(\ref{eq4}) can be formulated equivalently as:
\begin{equation}\label{prob_eq1}
\begin{aligned}
\min\limits_{(u,z)\in L^2(\mathcal{O})\times L^2(\mathcal{O})}\quad &F(u)+I_{\mathcal{C}}(z)\\
\mbox{s.t.}\quad &u=z
\end{aligned}
\end{equation}
where $F(u)$ denotes the objective functional of problem (\ref{eq3})-(\ref{eq4}) and $I_{\mathcal{C}}$ is the indicator function of admissable set $\mathcal{C}$. Because of the separable structure of problem (\ref{prob_eq1}), ADMM type methods \cite{MR3170605} can be obviously used to solve it.

Each iteration of ADMM includes the minimization of an unconstrained linear parabolic optimal control problem (denotes as u-subproblem) and projection onto the admissible control set (denotes as z-subproblem). Indeed, when we apply ADMM to solve this problem, these two types of constraints (\ref{eq2}), (\ref{eq4}) are treated separately. The z-subproblem is easy since it has the closed-form solution. But the u-subproblem can only be solved iteratively by some certain numerical method. For example, as studied in \cite{MR2404764}, we can use conjugate gradient method to solve it. Clearly, solving the u-subproblem dominates the computation amount of each iteration. Furthermore, the total dimension of the unconstrained parabolic optimal control problem could be very large after the discretization of space and time. Thus, it is time-consuming to obtain a high-precision solution of the u-subproblem at each iteration.

Recently, Song, eta., \cite{song2020implementation} have proposed an easily implementable and appropriately accurate inexactness criterion for solving u-subproblem at each ADMM iteration. Their method is called inexact ADMM which can be regarded as an improved version of ADMM method for solving problem (\ref{eq3})-(\ref{eq4}). Besides the convergence of inexact ADMM is proved. They also show efficiency of their method comparing with ADMM through some concrete examples. 
	
   The obtained u-subproblem at each iteration can be solved inexactly that means each iteration may be easier, but the convergence rate of ADMM type method in both ergodic and non-ergodic sense is order of $O(1/K)$ where $K$ denotes the iteration counter \cite{MR2914282,MR3347463}. Since slow convergence rate results in more iterations, inexact ADMM method sometimes may be not an efficient method. This is also the main defect of ADMM type methods.
      
   Notice that the indicator function of the additional constraint on control variable $u\in\mathcal{C}$ arises in the optimality condition of problem (\ref{eq3})-(\ref{eq4}). This fact motivates one to consider the semi-smooth Newton (SSN) methods because the indicator function is nonsmooth. Besides, the SSN type methods have been widely studied for elliptic optimal control problems with control constraints(e.g. \cite{MR3504560,MR2955288}). Naturally, we should consider whether they can be directly extended to solve the problem (\ref{eq3})-(\ref{eq4}).  
                  
      The SSN type methods' common feature is that the semi-smooth Newton direction is constructed by using a generalized Jacobian matrix in the sense of Clarke \cite{MR1058436} and then each iteration is expressed in terms of certain active set strategy which identifies the active and inactive indices with respect to box constraints, see \cite{MR3504560,MR1972219} for more detail. In \cite{MR1691937} a special semi-smooth Newton method with the active set strategy, called the primal-dual active set (PDAS) method is introduced for solving control constrained elliptic optimal control problems. The convergence result of PDAS approach can be founded in \cite{MR1951024} and this method can also be extended to solve parabolic boundary optimal control problems. It is proved in \cite{MR1972217} that SSN type methods possess a local superlinear convergence rate and the solution own high-precision as long as we have a good initial guess. 
	 
	 Although SSN type methods possess rigorous theoretical results, directly extending SSN type methods to solve our problem seems difficult to implement because of the following problems. First of all, the dimensionality of the resulting Newton system restricts implementation. For example, the case $n=2$, we set the mesh sizes of time and space discretization as $1/100$, then the dimensionality of Newton system obtained at each iteration is order of $O(10^{6})$ and $O(10^{8})$ for the case $n=3$. Secondly, SSN type methods require us to solve Newton system exactly otherwise convergence and convergence rate results cannot achieve in numerical computation. What's more, the system to
	 be solved at each iteration is large and ill-conditioned linear equations which mean preconditioner's design is required. 
	 When we apply SSN type methods, the box constraints (\ref{eq4}) on control variable are forced to be considered together with linear parabolic PDE constraint (\ref{eq2}) simultaneously. Hence the varying active set results  Newton equation of different structure at each iteration, which forces us to adjust the preconditioner at each step. Thus solving the obtained Newton system exactly at each iteration is hard to implement from both computational load and computational amount perspective.
	 
	  The main defect of ADMM type methods is the slow convergence rate, and the SSN type methods' implementation is mainly restricted by the high dimensionality of discretized problem. Our desire is to design some implementable and more efficient numerical schemes to solve this problem.
      
      We notice that no matter ADMM type methods or SSN type methods, their design is mainly based on how to solve the primal problem, though ADMM can be interpreted from the way of solving dual problem \cite{davis2016convergence} and SSN can be regarded as a primal-dual method \cite{MR1951024}. Hence we consider whether we can design algorithms by the way of solving the dual problem. Note that Burachik, eta., \cite{MR3199411} have studied the Fenchel dual problem of the control constrained optimal control problem. The dynamic system of their considered problem is described by linear ODEs. They showed that the solution of primal problem can be obtained by solving the dual problem. Besides Christian Clason, eta., \cite{MR2775195} have considered the dual problem of some specific unconstrained elliptic optimal control problems. These studies also motivate us to consider the dual problem of (\ref{eq3})-(\ref{eq4}).
      
      In the present paper, we derive the specifically dual problem for control bounded linear parabolic optimal control problem (\ref{eq3})-(\ref{eq4}), using Fenchel duality scheme. The dual problem is only constrained by linear parabolic PDE and the objective functional is first-order differentiable that means it is an unconstrained smooth optimal control problem. Besides, the control constraints in the primal problem are embedded in the dual problem's objective functional. The main advantage of the dual problem compared to primal problem is the vanished box constraints. This fact help us design some implementable numerical methods. 
    
    The rest of this paper is organized as follows. In Section 2, we first recall Fenchel duality relevant to our work then derive the dual problem for (\ref{eq3})-(\ref{eq4}).
     We also prove the existence and uniqueness of solution to the dual problem and derive the associated optimality conditions. Section 3 and Section 4 are concerned with numerical algorithm for solving the dual problem obtained by Section 2. In these two section, we discuss how to design first-order and second-order algorithms respectively. Besides, the numerical discretization of the dual problem by finite difference and finite element method is also discussed. Some preliminary numerical results for the algorithm designed by us are reported in Section 5 to validate the efficiency of our proposed numerical methods. Finally, we make concluding remarks and identify some future work in Section 6. 
	
	\textbf{Remark.} For convenience, in the rest of this paper we assume that the initial value $y_0=0$ in (\ref{eq2}). This assumption makes the solution operator w.r.t. pde is linear. For general case ($y_0\neq 0$), we can consider the principle of superposition of solutions and convert problem into the case $y_0=0$.
	\section{The dual problem}
	In this section, firstly we present some notations and known results corresponding to Fenchel duality that will be used in the later analysis. Then we derive the dual problem specifically, prove the existence and uniqueness of optimal control for dual problem. Finally the associated first-order optimality conditions are derived. 
	
	\subsection{Preliminaries}
	Here we briefly recall Fenchel duality, complete discussion can be found in \cite{MR1727362,MR1451876}. Let $V$ and $Y$ be Banach spaces with topological dual space $V^{*}$ and $Y^{*}$, respectively, and let $\Lambda: V \rightarrow Y$ be a continuous linear operator. Furthermore, $\Lambda^*$ represents the adjoint operator of $\Lambda$, $\Lambda^*: Y^* \rightarrow V^*$. 
	The following Theorem is called Fenchel duality theorem \cite{MR1727362}.
	\begin{theorem}
		Let $\mathcal{F}: V \rightarrow \overline{\mathbb{R}}$, $\mathcal{G}: Y \rightarrow \overline{\mathbb{R}}$ be convex lower semicontinuous functionals which are not identically equal $\infty$ and there exists some $v_{0} \in V$ such that $\mathcal{F}\left(v_{0}\right)<\infty, \mathcal{G}\left(\Lambda v_{0}\right)<\infty$, and $\mathcal{G}$ is continuous at $\Lambda v_{0}$, then there holds
	\begin{equation}\label{eq1}
	\inf _{ v \in V} \mathcal{F}(v)+\mathcal{G}(\Lambda v)=\sup _{q \in Y^{*}}-\mathcal{F}^{*}\left(\Lambda^{*} q\right)-\mathcal{G}^{*}(-q)
	\end{equation}
	 Furthermore, the equality in (\ref{eq1}) is attained at $\left(v^{*}, q^{*}\right)$ if and only if
	\begin{equation}
	\left\{\begin{array}{c}
	\Lambda^{*} q^{*} \in \partial \mathcal{F}\left(v^{*}\right) \\
	-q^{*} \in \partial \mathcal{G}\left(\Lambda v^{*}\right)
	\end{array}\right.
	\end{equation}
	\end{theorem}

Here $\mathcal{F}^{*}: V^{*} \rightarrow\overline{\mathbb{R}}$ denotes the Fenchel conjugate functional of $\mathcal{F}$ defined by
$$
\mathcal{F}^{*}(q)=\sup\limits_{v \in V}\langle q,v\rangle_{V^{*}, V}-\mathcal{F}(v)
$$
where $\langle q,v\rangle_{V^{*}, V}:=q(v)$

There holds the following important fact:
$$
\mathcal{F}^{*}(q)=\langle q,v\rangle_{V^{*}, V}-\mathcal{F}(v) \quad \text { if and only if } \quad q \in \partial \mathcal{F}(v)
$$
Here, $\partial \mathcal{F}$ denotes the subdifferential of the convex functional $\mathcal{F}$, which reduces to the Gâteaux-derivative if it is Gâteaux differentiable.

	\subsection{Derivation of dual problem}
	In this subsection we derive the dual problem associated with problem (\ref{eq3})-(\ref{eq4}). The main theoretical tool is Fenchel duality discussed in the previous subsection.
	
	We introduce the linear operator  $S: L^2(\mathcal{O})\rightarrow L^2(Q)$ associated with state equation (\ref{eq2}), and it is defined as
	\begin{equation*}
	S(u):=y
	\end{equation*}
	It is shown in \cite{MR0271512} that $S$ is continuous, compact and inversable. We denote $S^*$ as the adjoint operator of $S$. 
	
	Then optimal control problem (\ref{eq3})-(\ref{eq4}) can be formulated equivalently as
	\begin{equation}
	\min\limits_{u} \frac{1}{2} \iint_{Q}\left|S(u)-y_{d}\right|^{2} d x d t+\frac{\gamma}{2} \iint_{\mathcal{O}}|u|^{2} dx dt+I_{\mathcal{C}}(u)
	\end{equation}
	where $I_{\mathcal{C}}(\cdot)$ denotes indicator function of the admissble set $\mathcal{C}$ that is,
	\begin{equation*}
	I_{\mathcal{C}}(z)=\left\{
	\begin{aligned}
	0, \quad&\text{if }z\in\mathcal{C}\\
	+\infty, \quad&\text{if }z\in L^2(Q)\setminus\mathcal{C}
	\end{aligned}
	\right.
	\end{equation*}
	
	Since functional space $L^2(\mathcal{O})$ and $L^2(Q)$ is reflexive, we can directly calculate the dual problem. In order to follow Fenchel duality scheme, we define the following functional:
	\begin{align}
	&\mathcal{F}:L^2(\mathcal{O})\rightarrow \overline{\mathbb{R}} \quad  \mathcal{F}(u)=\frac{\gamma}{2} \iint_{\mathcal{O}}|u|^{2} dx dt+I_{\mathcal{C}}(u)\\
	&\mathcal{G}:L^2(Q)\rightarrow\overline{\mathbb{R}} \quad \mathcal{G}(y)=\frac{1}{2} \iint_{Q}\left|y-y_{d}\right|^{2} d x d t
	\end{align}
	
    The Fenchel conjugate functional of $\mathcal{F}$ and $\mathcal{G}$ are given by:
	\begin{align}
&\mathcal{F}^*:L^2(\mathcal{O})\rightarrow \overline{\mathbb{R}} \quad	\mathcal{F}^*(p)=\langle p,\operatorname{Pr}_{\mathcal{C}}\left(\frac{p}{\gamma}\right)\rangle_{L^2(\mathcal{O})}-\frac{\gamma}{2}\|\operatorname{Pr}_{\mathcal{C}}\left(\frac{p}{\gamma}\right)\|^2_{L^2(\mathcal{O})}\\
 &\mathcal{G}^*:L^2(Q)\rightarrow \overline{\mathbb{R}} \quad   \mathcal{G}^*(q)=\langle q,y_d\rangle_{L^2(Q)}+\frac{1}{2}\|q\|^2_{L^2(Q)}
	\end{align}
	where $\operatorname{Pr}_{\mathcal{C}}(\cdot)$ denotes the projection onto the admissible set $\mathcal{C}$, mathematically,
	\begin{equation*}
	\operatorname{Pr}_{\mathcal{C}}(v(x,t))=\left\{
	\begin{aligned}
	&v(x,t), \quad&\text{if }v(x,t)\in\mathcal{C}\\
	&a, \quad&\text{if }v(x,t)<a\\
	&b, \quad&\text{if }v(x,t)>b\\
	\end{aligned}
	\right.
	\end{equation*}
		
	Since $\mathcal{F}$ and $\mathcal{G}$ are convex and lower semi-continuous, $S$ is also a continuous linear operator, the dual problem conceptually can be formulated by $\max\limits_{q \in L^2(Q)}-\mathcal{F}^{*}\left(S^{*} (q)\right)-\mathcal{G}^{*}(-q)$. We formulate dual problem  equivalently as following:
	\begin{equation}\label{dual}
	\begin{aligned}
	\min\limits_{\bs q\in L^2(Q)}\quad &\langle p,\operatorname{Pr}_{\mathcal{C}}\left(\frac{p}{\gamma}\right)\rangle_{L^2(\mathcal{O})}-\frac{\gamma}{2}\|\operatorname{Pr}_{\mathcal{C}}\left(\frac{p}{\gamma}\right)\|^2_{L^2(\mathcal{O})}+\langle -q,y_d\rangle_{L^2(Q)}+\frac{1}{2}\|q\|^2_{L^2(Q)}\\
	\mbox{s.t.}\quad &p=S^{*}(q)
	\end{aligned}
	\end{equation}
	
	In order to specify dual problem, we must express constraint $p=S^*(q)$ explicitly. Since solution operator $S$ is inversable, the constraint $p=S^{*}(q)$ is equivalent to $q=S^{-*}(p)$ where $S^{-*}$ denotes the inverse operator of operator $S^*$.
	
	By the definition of adjoint operator there holds $\langle S^{-1}y, p\rangle_{L^2(\mathcal{O})}=\langle y,S^{-*}p\rangle_{L^2(Q)}$. Then taking advantage of PDE constraint (\ref{eq2}), we can derive that
	
	\begin{equation}\label{eq5}
	\begin{aligned}
	\langle S^{-1}y, p\rangle_{L^2(\mathcal{O})}&=\langle u\chi_{\mathcal{O}}, p \rangle_{L^2(Q)}\\
	&=\langle\frac{\partial y}{\partial t}-\nu\Delta y+a_{0} y,p\rangle_{L^2(Q)}
	\end{aligned}
	\end{equation}
	Intergration by parts in time and application Green's formula in space finally yields the following equation:
	\begin{equation*}
	\langle\frac{\partial y}{\partial t}-\nu\Delta y+a_{0} y,p\rangle_{L^2(Q)}=\iint_{\Omega} y(x,T)\cdot p(x,T)dx-\langle y,\frac{\partial p}{\partial t}+\nu\Delta p-a_0p\rangle_{L^2(Q)}-\nu\iint_{\Gamma\times(0,T)}(\frac{\partial y}{\partial n}p-\frac{\partial p}{\partial n}y)dxdt
	\end{equation*} 
	We set teriminal condition $p(x,T)=0$ and $p=0$ on $\Gamma\times(0,T)$, combining with (\ref{eq5}) there holds:
	\begin{equation*}
	\langle S^{-1}y, p\rangle_{L^2(\mathcal{O})}=-\langle y,\frac{\partial p}{\partial t}+\nu\Delta p-a_0p\rangle_{L^2(Q)}
	\end{equation*}
	By the definition of adjoint operator we conclude 
	\begin{equation*}
	S^{-*}p=-(\frac{\partial p}{\partial t}+\nu\Delta p-a_0p)
	\end{equation*}	
	
	For convenience, we denote the objective functional of dual problem as $J(q)$. The dual problem (\ref{dual}) can be formulated explicitly:
	\begin{equation}\label{dual_prob}
\min\limits_{q} J(q):=\langle p,\operatorname{Pr}_{\mathcal{C}}\left(\frac{p}{\gamma}\right)\rangle_{L^2(\mathcal{O})}-\frac{\gamma}{2}\|\operatorname{Pr}_{\mathcal{C}}\left(\frac{p}{\gamma}\right)\|^2_{L^2(\mathcal{O})}+\langle -q,y_d\rangle_{L^2(Q)}+\frac{1}{2}\|q\|^2_{L^2(Q)}
	\end{equation}
	subject to the state equation:
	\begin{equation}\label{dual_prob_con}
	\left\{
	\begin{aligned}
	&\frac{\partial p}{\partial t}+\nu\Delta p-a_0p=-q\quad&\text{in}\quad\Omega\times(0,T)\\
	&p=0\quad&\text{on}\quad\Gamma\times(0,T)\\
	&p(T)=0
	\end{aligned}
	\right.
	\end{equation}
	
We can treat variable $q$ as control variable and variable $p$ as state variable. Thus dual problem (\ref{dual_prob})-(\ref{dual_prob_con}) is an unconstrained parabolic optimal control problem.

The differentiability of objective functional is an important
property for the numerical optimization method design. Thus in the final of this subsection, we prove the fact that the dual problem's objective functional $J(q)$ is continuously differentiable.
 
	\begin{theorem}\label{smooth}
		The objective functional $J(\cdot)$ is Gâteaux differentiable corresponding to state variable $p$ and control variable $q$ respectively.
	\end{theorem}
	
	{\emph{Proof.}} 
	It is obviously that objective functional $J(\cdot)$ is differentiable corresponding to variable $q$. Thus we just need to show that it is Gâteaux differentiable associated with variable $p$.
	
	For convenience, we introduce function $\theta(x)$ where $\theta: \mathbb{R}\rightarrow \mathbb{R}$ is defined by:
	\begin{equation}\label{theta}
	\theta(x)=\left\{\begin{aligned}
	&\frac{x^2}{2\gamma}\quad &\text{if}\quad a\leq\frac{x}{\gamma}\leq b\\
	&a\cdot x-\frac{\gamma}{2}a^2\quad&\text{if}\quad \frac{x}{\gamma}<a\\
	&b\cdot x-\frac{\gamma}{2}b^2\quad&\text{if}\quad \frac{x}{\gamma}>b 
	\end{aligned}
	\right.
	\end{equation}
		
	Thus, there must hold the following equation
	\begin{equation*}
	\langle p,\operatorname{Pr}_{\mathcal{C}}\left(\frac{p}{\gamma}\right)\rangle_{L^2(\mathcal{O})}-\frac{\gamma}{2}\|\operatorname{Pr}_{\mathcal{C}}\left(\frac{p}{\gamma}\right)\|^2_{L^2(\mathcal{O})}=\iint_{\mathcal{O}}\theta(p(x,t))dxdt
	\end{equation*}
	
	Besides we notice that $\theta: \mathbb{R}\rightarrow\mathbb{R}$ is continuously differentiable. Combining with the definition of Gâteaux differentiable there holds that $J$ is also Gâteaux differentiable associated with variable $p$. Thus we complete the proof. $\hfill\qedsymbol$  
	
	\textbf{Remark.} The objective functional of dual problem is Gâteaux differentiable while box constraints (\ref{eq4}) in primal problem are embedded in the dual objective functional (\ref{dual_prob}). Unlike the constrained control variable in primal problem, the dual problem's control variable is unconstrained.
	\subsection{Existence of optimal control for the dual problem}
	We prove in this subsection the existence of optimal control for the dual problem. By this fact, we can further discuss the relationship between the solution of primal problem and the solution of dual problem (saddle point property). That means we can get the solution of primal problem (\ref{eq3})-(\ref{eq4}) by solving its dual problem (\ref{dual_prob})-(\ref{dual_prob_con}).
	\begin{theorem}
		 There exists a unique optimal control $\bar{q}\in L^{2}\left(Q\right)$ such that $J(\bar{q}) \leq J(q), \forall q\in L^{2}\left(Q\right)$
	\end{theorem}
	 	
	{\emph{Proof.}} Firstly, we observe that 
	\begin{equation*}
	 \langle p,\operatorname{Pr}_{\mathcal{C}}\left(\frac{p}{\gamma}\right)\rangle_{L^2(\mathcal{O})}-\frac{\gamma}{2}\|\operatorname{Pr}_{\mathcal{C}}\left(\frac{p}{\gamma}\right)\|^2_{L^2(\mathcal{O})}=\iint_{\mathcal{O}}\theta(p(x,t))dxdt\geq 0
	\end{equation*}
	where $\theta(\cdot)$ is the function defined by (\ref{theta}). The above inequality holds because $\theta(\cdot)\geq 0$. Then we obtain the following inequality:
	\begin{equation}\label{ineq1}
	\begin{aligned}
	J(q) &\geq \langle -q,y_d\rangle_{L^2(Q)}+\frac{1}{2}\|q\|^2_{L^2(Q)}\\ &=\frac{1}{2}\|q-y_d\|^2_{L^2(Q)}-\frac{1}{2}\|y_d\|_{L^2(Q)}^2\geq -\frac{1}{2}\|y_d\|_{L^2(Q)}^2,\quad \forall q\in L^{2}\left(Q\right)
	\end{aligned}
	\end{equation} 
	
	Thus the infimum of $J(q)$ exists and there must exist a sequence $\left\{q_{n}\right\} \subsetneq L^{2}\left(Q\right)$ such that
	$$
	\lim_{n \rightarrow \infty} J\left(q_{n}\right)=\inf_{q\in L^{2}\left(Q\right)} J(q)
	$$
	Combining with (\ref{ineq1}), precisely
	\begin{equation*}
	 \frac{1}{2}\left\|q_{n}-y_d\right\|_{L^2(Q)}^{2} \leq J\left(q_{n}\right)+\frac{1}{2}\|y_d\|_{L^2(Q)}^2
	 \end{equation*}
	  implies that $\left\{q_{n}\right\}$ is bounded in $L^{2}\left(Q\right)$. 
	  
	  Since $L^{2}\left(Q\right)$ is reflexive Hilbert space, there exists a subsequence of $\{ q_n\}$, still denoted by $\left\{q_{n}\right\}$, that converges weakly to $\bar{q}$ in $L^{2}\left(Q\right)$. 
	  
	  Because $J$ is convex and continuous, it must be weakly lower semi-continuous. Thus there holds:
	\begin{equation*}
	J(\bar{q}) \leq \liminf_{n \rightarrow \infty} J\left(q_{n}\right)=\inf_{q\in L^{2}\left(Q\right)} J(q)
	\end{equation*}
	
	We must have $\bar{q}$ is an optimal control for dual problem.

		Besides the uniqueness of optimal control $\bar{q}$ can be easily guaranteed because objective functional $J$ is strictly convex corresponding to control variable $q$. Thus we complete the proof.
	$\hfill\qedsymbol$\\
	
Since $\mathcal{F}$ and $\mathcal{G}$ are convex and lower semi-continuous functionals, solution operator $S$ is a continuous linear operator, the Fenchel duality theorem holds. The existence of solution to dual problem and primal problem guarantees that the optimal value is attainable. We set primal problem solution pair as $(\bar{y},\bar{u})$ and dual solution pair as $(\bar{p},\bar{q})$. Then we would like to find the relationship between $(\bar{y},\bar{u})$ and $(\bar{p},\bar{q})$. 

According to Fenchel duality theorem, the following general equations are satisfied:
	\begin{equation}
	\left\{\begin{array}{c}
	S^{*}(\bar{q}) \in \partial \mathcal{F}\left(\bar{u}\right) \\
	-\bar{q} \in \partial \mathcal{G}\left(S(\bar{u})\right)
	\end{array}\right.
	\end{equation}
	
   Since $\bar{p}=S^*(\bar{q})$, there holds
   \begin{equation*}
   \begin{aligned}
  &0\in\partial(\mathcal{F}(\bar{u})-\langle \bar{p},\bar{u}\rangle_{L^2(\mathcal{O})})\\
  \iff&u^*=\mathop{\arg\min}\limits_{u\in L^2(\mathcal{O})}\{\frac{\gamma}{2}\|u\|_{L^2(\mathcal{O})}^2-\langle \bar{p},u\rangle_{L^2(\mathcal{O})}+I_{\mathcal{C}}(u)\}\\
  \iff&\bar{u}=\operatorname{Pr}_{\mathcal{C}}\left(\frac{\bar{p}}{\gamma}\right)
   \end{aligned} 
   \end{equation*}
   
   Similarly one can also derive that $\bar{y}=y_d-\bar{q}$, thus the solution of primal problem can be obtained by the solution of dual problem:
   \begin{equation}\label{rel_pri_dual}
   \left\{\begin{aligned}
   \bar{u}&=\operatorname{Pr}_{\mathcal{C}}\left(\frac{\bar{p}}{\gamma}\right)\\
   \bar{y}&=y_d-\bar{q}
   \end{aligned}\right.
   \end{equation}
   
	\subsection{First-order optimality conditions for the dual problem}
	Let $DJ(q)$ be the first-order differential of $J$ at $q$ and $ \bar{q}$ the unique optimal control for the dual problem. Then it must hold the first-order optimality condition at optimal solution
	\begin{equation*}
	DJ(\bar{q})=0
	\end{equation*}
	In the rest of this subsection, we discuss the computation of gradient $DJ(q)$ that plays an important role in subsequent section.
	
	To compute $DJ(q)$, we employ a formal perturbation analysis as described in \cite{MR2404764}. Let $\delta q\in L^2(Q)$ be a peturbation of some $q\in L^2(Q)$, there holds
	\begin{equation}\label{derivative}
	\delta J(q)=\iint_{Q}DJ(q)\delta qdxdt
	\end{equation}
    and also holds
	\begin{equation}\label{perb_a_eq1}
	\delta J(q)=\langle\operatorname{Pr}_{\mathcal{C}}\left(\frac{p}{\gamma}\right),
	\delta p\rangle_{L^2(\mathcal{O})}-\langle y_d-q,\delta q\rangle_{L^2(Q)}
	\end{equation}
	in which $\delta p$ is the solution of
	\begin{equation}\label{perb_a_eq2}
	\left\{
	\begin{aligned}
	&\frac{\partial\delta p}{\partial t}+\nu\Delta \delta p-a_0\delta p=-\delta q\quad&\text{in}\quad\Omega\times(0,T)\\
	&\delta p=0\quad&\text{on}\quad\Gamma\times(0,T)\\
	&\delta p(T)=0
	\end{aligned}
	\right.
	\end{equation}
	
	Consider function $z$ defined over $Q$ and function $z$ is differentiable corresponding to variable $x$ and $t$. We multiple both sides of the first equation in (\ref{perb_a_eq2}) by function $z$ and integrate over $Q$. Then integration by parts in time and application of Green's formula in space finally yields
	\begin{equation}\label{perb_a_eq3}
	\begin{aligned}
		&\int_{\Omega} \delta p(T)z(T) d x-\int_{\Omega}\delta p(0) z(0) d x+\iint_{Q}\left[-\frac{\partial z}{\partial t}+\nu \Delta z-a_{0} z\right]\delta p d x d t \\
		&+\nu \iint_{\Gamma \times(0, T)}\left(\frac{\partial z}{\partial n}\delta p-\frac{\partial\delta p}{\partial n} z\right) d x d t=-\iint_{Q}y \delta q d x d t .
		\end{aligned}
	\end{equation}
	
	Let us assume $z$ is the solution to the adjoint parabolic equation
	\begin{equation}\label{perb_a_eq4}
	\left\{
	\begin{aligned}
	&\frac{\partial z}{\partial t}-\nu \Delta z+a_{0} z=\operatorname{Pr}_{\mathcal{C}}\left(\frac{p}{\gamma}\right)\chi_{\mathcal{O}}\quad&\text{in}\quad\Omega\times(0,T)\\
	& z=0\quad&\text{on}\quad\Gamma\times(0,T)\\
	& z(0)=0
	\end{aligned}
	\right.
	\end{equation}
	By equations (\ref{perb_a_eq1})-(\ref{perb_a_eq4}), there holds
\begin{equation*}
	\delta J(q)=\langle z-y_d+q,\delta q\rangle_{L^2(Q)}
	\end{equation*}
	together with (\ref{derivative}) we obtain
	\begin{equation}
	DJ(q)=z-y_d+q
	\end{equation}
	
	Thus, the first-order optimality conditions for the dual problem can be summarized as follows:
	\begin{theorem}
		Let $\bar{q}$ be the unique solution of dual problem (\ref{dual_prob})-(\ref{dual_prob_con}). Then it is characterized by the following optimality conditions:
		\begin{equation}\label{opt_sys1}
		DJ(\bar{q})=\bar{z}-y_d+\bar{q}=0
		\end{equation}
		where $\bar{z}$ is obtained from the successive solution of the following two parabolic PDE:
		\begin{equation}\label{opt_sys2}
		\begin{aligned}
		\text{(state equation)}&\left\{
		\begin{aligned}
		&\frac{\partial\bar{p}}{\partial t}+\nu\Delta \bar{p}-a_0\bar{p}=-\bar{q}\quad&\text{in}\quad\Omega\times(0,T)\\
		&\bar{p}=0\quad&\text{on}\quad\Gamma\times(0,T)\\
		&\bar{p}(T)=0
		\end{aligned}
		\right.\\
		\text{(adjoint equation)}&\left\{
		\begin{aligned}
		&\frac{\partial\bar{z}}{\partial t}-\nu \Delta\bar{z}+a_{0} \bar{z}=\operatorname{Pr}_{\mathcal{C}}\left(\frac{\bar{p}}{\gamma}\right)\chi_{\mathcal{O}}\quad&\text{in}\quad\Omega\times(0,T)\\
		&\bar{z}=0\quad&\text{on}\quad\Gamma\times(0,T)\\
		&\bar{z}(0)=0
		\end{aligned}
		\right.
		\end{aligned}
		\end{equation}
	\end{theorem}
	
	\section{First-Order Algorithm Design}
	In this section, we discuss how to numerically solve the dual problem(\ref{dual_prob})-(\ref{dual_prob_con}). Since dual problem is an unconstrained smooth optimal control problem, any first-order optimization method can be applied to solve it. Considering the problem is of large-scale after fully discretization, here we discuss the application of FR conjugate gradient method \cite{MR2244940} for problem(\ref{dual_prob}) and elaborate on the computation of the gradient and stepsize at each iteration. Finally an easily implementable numerical algorithm is obtained. 
	\subsection{FR-CG Framework for dual problem}
	Conceptually, the following FR-CG algorithm can be implemented to solve the dual problem.
	\begin{enumerate}
	\item[\textbf{(a)}] Give an initial guess $q^0$ and stopping tolerence constant $tol$. Set iteration number $k=0$. 
	\item[\textbf{(b)}] Compute the gradient of objective functional at $q^0$ denoted as $g^0=DJ(q^0)$ by solving the state equation (\ref{eq2}) and adjoint equation (\ref{perb_a_eq4}) corresponding to $q^0$. If $\|g^0\|_{L^2(Q)}<tol$, then set $\bar{q}=q^0$ as solution; otherwise set $d^0=-g^0$. For $k\geq 0$, we compute $q^{k+1}$, $g^{k+1}$ and $d^{k+1}$ as follows:
	\item[\textbf{(c)}] Choose the stepsize $\rho_k$ by solving the following optimization problem which is called exact line-search:
    \begin{equation}\label{exact_line}
    \left\{
    \begin{aligned}
    &\rho_k\in\mathbb{R}\\
    &J(q^k+\rho_kd^k)\leq J(q^k+\rho d^k) \quad \forall \rho\in\mathbb{R}
    \end{aligned}
    \right.
    \end{equation}

	\item[\textbf{(d)}] Update $q^{k+1}$ and $g^{k+1}$ respectively by
	\begin{equation*}
	\begin{aligned}
	q^{k+1}=q^k+\rho_kd^k\\
	g^{k+1}=DJ(q^{k+1})
	\end{aligned}
	\end{equation*}
	
	If $\|g^{k+1}\|_{L^2(Q)}<tol$, take $\bar{q}=q^{k+1}$ as the solution otherwise compute 
	\begin{equation*}
	\beta_k=\frac{\|g^{k+1}\|^2_{L^2(Q)}}{\|g^k\|^2_{L^2(Q)}}
	\end{equation*}
	and then update $d^{k+1}$ by
	\begin{equation*}
	d^{k+1}=-g^{k+1}+\beta_kd^k
	\end{equation*}
	
	Set $k=k+1$ and return to \textbf{(c)}.
	\end{enumerate}
	
	The above iterative framework looks quite simple, but it is formulated in functional spaces. In order to obtain an implementable framework, we need discretize the obtained dual problem. Besides the gradient computation of the discretized optimization problem and the stepsize choice in each iteration is important for FR-CG method. We shall discuss how to approach these two issues in the rest of this section. 
	\subsection{Computation of stepsize $\rho_k$}
	An important issue of the FR-CG method described by (\textbf{a})-(\textbf{d}) is the computation of the stepsize $\rho_k$ at each step. Using (\ref{exact_line}) to determine stepsize is numerically expensive and hard to implement because the objective functional is nonlinear and nonquadratic that means the explict formular for stepsize cannot be derived. If we want to get the stepsize by formular (\ref{exact_line}), we can only advocate iterative method for solving the following equation:  
	\begin{equation*}
	\langle DJ(q^k+\rho^k\cdot d^k),d^k \rangle=0
	\end{equation*}
	
	However one should notice that at each iteration for solving the above equation the gradient should be evaluated which means two parabolic equations required to be solved. Hence, the implementation of exact line-search for stepsize computation is numerically expensive.
	
	The high computational load for solving (\ref{exact_line}) motivates us to implement certain stepsize rule to find an 
	appropriate stepsize $\rho_k$. Here we advocate the following inexact line-search rule that strategy is called Armijo condition \cite{MR2244940}.
	\begin{equation}\label{inexact_line}
	J(q^k+\rho d^k)\leq J(q^k)+c\cdot\langle DJ(q^k),d^k\rangle_{L^2(Q)}\cdot \rho
	\end{equation}
	where $c$ is a given constant. 
	
	If we choose (\ref{inexact_line}) to determine stepsize, the objective functional value is required to be evalulated repeatedly. For a given $d^k\in L^2(Q)$, we find that the state $p=S^*(q^k+\rho d^k)$ in the objective functional $J(q^k+\rho d^k)$ can be computed by the following equation that is due to operator $S^*$ is linear:
	\begin{equation*}
	p=S^*(q^k+\rho d^k)=S^*(q^k)+\rho S^*(d^k)
	\end{equation*}
	
   	Furthermore, that means for any $\rho\in\mathbb{R}$ the evaluation of functional value $J(q^k+\rho d^k)$ just need solve one more parabolic equation w.r.t. $d^k$. Hence, repeatedly evaluating functional value is implementable.

	For comparsion, the stepsize determine procedure by Armijo condition just need solve one more parabolic equation, while the iteration method for formular (\ref{exact_line}) need at least solve four parabolic equation since we cannot guarantee initial guess to be good enough.   
	\subsection{Problem discretization and gradient computation}
	In this subsection, we first discuss the numerical discretization for our dual problem(\ref{dual_prob})-(\ref{dual_prob_con}), thus we can obtain a large-scale finite-dimensional optimization problem. We employ the backward-Euler finite difference method for time discretization and piecewise linear element method for space discretization. The computation of gradient w.r.t. discretized problem is also considered. Finally an easily implementable FR-CG method for the fully discrete dual problem is obtained.

	Firstly, we discuss time discretization technique for dual problem.
	
	We define the time step $\Delta t$ by $\Delta t=\frac{T}{N}$, with $N$ a positive integer. Thus, we approximate the admissible control space  $L^2(Q)$ by $[L^2(\Omega)]^{N}$; and equip $[L^2(\Omega)]^{N}$ with the following inner product
	$$
	\langle v, w\rangle_{\Delta t}=\Delta t \sum_{n=0}^{N-1} \int_{\Omega} v_{n} \cdot w_{n} d x, \quad \forall v=\left\{v_{n}\right\}_{n=0}^{N-1}, w=\left\{w_{n}\right\}_{n=0}^{N-1} \in [L^2(\Omega)]^{N}
	$$
	and the norm
	$$
	\|v\|_{\Delta t}=\left(\Delta t \sum_{n=0}^{N-1} \int_{\Omega}\left|v_{n}\right|^{2} d x\right)^{\frac{1}{2}}, \quad \forall v=\left\{v_{n}\right\}_{n=0}^{N-1} \in [L^2(\Omega)]^{N}
	$$

    Then the original problem(\ref{dual_prob})-(\ref{dual_prob_con}) is approximated by the following semi-discrete optimal control problem.
	\begin{equation}\label{discrete_prob}
	\left\{\begin{aligned}
	&\bar{q}^{\Delta t}\in [L^2(\Omega)]^{N}\\
	&J^{\Delta t}\left(\bar{q}^{\Delta t}\right) \leq J^{\Delta t}(q^{\Delta t}),\quad \forall q^{\Delta t}=\left\{q_{n}\right\}_{n=0}^{N-1} \in [L^2(\Omega)]^{N}
	\end{aligned}\right.
	\end{equation}
	where the cost functional $J^{\Delta t}(q^{\Delta t})$ is defined by
	$$
	J^{\Delta t}(q^{\Delta t})=\langle p^{\Delta t},\operatorname{Pr}_{\mathcal{C}}\left(\frac{p^{\Delta t}}{\gamma}\right)\rangle_{\Delta t}-\frac{\gamma}{2}\|\operatorname{Pr}_{\mathcal{C}}\left(\frac{p^{\Delta t}}{\gamma}\right)\|^2_{\Delta t}+\langle -q^{\Delta t}, y_d^{\Delta t}\rangle_{\Delta t}+\frac{1}{2}\|q^{\Delta t}\|^2_{\Delta t}
	$$
	with $y_d^{\Delta t}=\left\{y_d^{n}\right\}_{n=0}^{N-1}$ and $y_d^{n}:=y_d(\cdot,n\Delta t)$, $p^{\Delta t}=\left\{p_{n}\right\}_{n=0}^{N-1}$ the solution of the following semi-discrete state equation: $p_{N}=0$ then for $n=0,1,\ldots,N-1$, with $p_{n+1}$ being known, we obtain $p_{n}$ from the solution of the following linear elliptic equation:
	$$
	\left\{\begin{aligned}
	\frac{p_{n}-p_{n+1}}{\Delta t}-\nu\Delta p_{n}+a_{0} p_{n}=q_{n} & \text { in }\Omega \\
	p_{n}=0 & \text { on }\Gamma
	\end{aligned}\right.
	$$
	
	Thus we should need to solve a simple elliptic problem to obtain $p_n$ from $p_{n+1}$. Besides, our scheme is first-order accurate and robust w.r.t. time variable. The existence of solution to semi-discretization problem (\ref{discrete_prob}) can be proved as done in dual problem.  
	
	Then we should discuss how to compute the gradient associated with problem (\ref{discrete_prob}).
	
	Let $\bar{q}^{\Delta t}$ be the solution to the semi-discrete problem, then it should satisfy the following necessary condition:
	\begin{equation*}
	DJ^{\Delta t}(\bar{q}^{\Delta t})=0
	\end{equation*}
   	
    Proceeding as in the continuous case, we can derive the gradient w.r.t. $q^{\Delta t}=\{q_n\}_{n=0}^{N-1}$ for discretized problem (\ref{discrete_prob}).
    \begin{equation*}
    DJ^{\Delta t}(q^{\Delta t})=\{z_n-y^n_d+q_n\}_{n=0}^{N-1}
    \end{equation*}

   Here $\{z_n\}_{n=0}^{N-1}$ is the solution of the following semi-discrete adjoint system:
   
   for $n=0$ we need solve the following linear ellptic equation, where $\chi_{O}$ represents characteristic function of the set $O$:
   \begin{equation*}
   \left\{
   \begin{aligned}
   \frac{z_0}{\Delta t}-\nu\Delta z_0+a_0z_0&=\operatorname{Pr}_{C}(\frac{p_0}{\gamma})\chi_{O} & \text { in }\Omega \\
   z_0&=0 & \text { on }\Gamma  
    \end{aligned}\right.
   \end{equation*}
   
   and for $n=1,\ldots,N-1$, with $z_{n-1}$ being known, we get $z_n$ by solving the following linear ellptic equation, where $\chi_{O}$ represents characteristic function of the set $O$:  
   \begin{equation*}
   \left\{
   \begin{aligned}
   \frac{z_{n}-z_{n-1}}{\Delta t}-\nu\Delta z_{n}+a_0z_{n}&=\operatorname{Pr}_{C}(\frac{p_{n}}{\gamma})\chi_{O} & \text { in } \Omega \\
   z_n&=0 & \text { on } \Gamma  
   \end{aligned}
   \right.
   \end{equation*}
   
   Next we shall discuss the space discretization and thus obtain a fully discrete problem. The gradient computation for fully discrete problem is also considered. For simplicity, we assume that $\Omega$ is a polygonal domain of $\mathbb{R}^2$ (since more complicated domain can be approximated by a family of such domain).
        
    Let $\mathcal{T}_h$ be a triangulation of $\Omega$ and $P_1$ the space of polynomial functions of two variables of degree less than one. We define finite element space $V_h$ and its subspace $V_{0h}$ by
    \begin{equation*}
    \begin{aligned}
    V_{h}&=\left\{\varphi_{h}\left|\varphi_{h} \in C^{0}(\bar{\Omega}) ; \varphi_{h}\right|_{\mathbb{T}} \in P_{1}, \forall \mathbb{T} \in \mathcal{T}_{h}\right\}\\
    V_{0 h}&=\left\{\varphi_{h}\left|\varphi_{h} \in V_{h}, \varphi_{h}\right|_{\Gamma}=0\right\}:=V_{h} \cap H_{0}^{1}(\Omega)
    \end{aligned}
    \end{equation*}
     
     Thus the semi-discretized control space $[L^2(Q)]^N$ is further approximated by $[V_h]^N$.
     
     The fully discrete optimal control problem that approximate the dual problem (\ref{dual_prob}) can be defined by:
     \begin{equation}\label{full_discrete_prob}
     \left\{\begin{aligned}
     &\bar{q}_h^{\Delta t}\in [V_h]^{N}\\
     &J_h^{\Delta t}\left(\bar{q}_h^{\Delta t}\right) \leq J^{\Delta t}_h(q_h^{\Delta t}),\quad \forall q_h^{\Delta t}=\left\{q_{n,h}\right\}_{n=1}^{N} \in [V_h]^{N}
     \end{aligned}\right.
     \end{equation}
     where the fully discrete cost functional $J_h^{\Delta t}$ is defined by
     \begin{equation}
     \begin{aligned}
     J^{\Delta t}_h(q_h^{\Delta t})=&\langle p_h,\operatorname{Pr}_{\mathcal{C}}\left(\frac{p_h}{\gamma}\right)\rangle_{\Delta t}-\frac{\gamma}{2}\|\operatorname{Pr}_{\mathcal{C}}\left(\frac{p_h}{\gamma}\right)\|^2_{\Delta t}+\langle -q_h, y_d^h\rangle_{\Delta t}+\frac{1}{2}\|q_h\|^2_{\Delta t}\\
     =&\Delta t\sum\limits_{n=0}\limits^{N-1}\int_{\Omega}p_{n,h}\cdot\operatorname{Pr}_{\mathcal{C}}\left(\frac{p_{n,h}}{\gamma}\right)dx-\frac{\gamma}{2}\Delta t\sum\limits_{n=0}\limits^{N-1}\int_{\Omega}\|\operatorname{Pr}_{\mathcal{C}}\left(\frac{p_{n,h}}{\gamma}\right)\|^2dx\\
     &-\Delta t\sum\limits_{n=0}\limits^{N-1}\int_{\Omega}q_{n,h}\cdot y_d^{n,h}dx+\frac{\Delta t}{2}\sum\limits_{n=0}\limits^{N-1}\int_{\Omega}\|q_{n,h}\|^2dx
     \end{aligned}
     \end{equation}
    with $\left\{ p_{n,h}\right\}_{n=0}^{N-1}$ the solution of the following fully discrete state equation: $p_{N,h}=0$ then for $n=0,1,\ldots,N-1$, with $p_{n+1,h}$ being known, we obtain $p_{n,h}$ from the solution of the following linear variational problem:
    \begin{equation}\label{opt_condition2}
    \left\{\begin{aligned}
    &p_{n,h}\in V_{0,h}\\
    &\int_{\Omega}\frac{p_{n,h}-p_{n+1,h}}{\Delta t}\varphi dx+\int_{\Omega}\nu\nabla p_{n,h}\cdot\nabla\varphi dx+\int_{\Omega}a_{0} p_{n,h}\varphi dx=\int_{\Omega}q_{n,h}\varphi dx\quad\text{for all\quad}\varphi\in V_{0,h} \\
    \end{aligned}\right.
    \end{equation}
    
    We can show that the first-order differential of $J_h^{\Delta t}$ at $q^{\Delta t}_h\in [V_h]^N$ is    
    \begin{equation}\label{opt_condition1}
    DJ^{\Delta t}_h(q^{\Delta t}_h)=\{z_{n,h}-y^{n,h}_d+q_{n,h}\}_{n=0}^{N-1}
    \end{equation}
    and $\{z_{n,h}\}_{n=0}^{N-1}$ is the solution of the following fully discrete adjoint system:
    
    for $n=0$ we need solve the following linear-variational problem
    \begin{equation}\label{opt_condition3}
    \left\{
    \begin{aligned}
    &z_{0,h}\in V_{0,h}\\
    &\int_{\Omega}\frac{z_{0,h}}{\Delta t}\varphi dx+\nu\int_{\Omega}\nabla z_{0,h}\cdot\nabla \varphi dx+\int_{\Omega}a_0z_{0,h}\varphi dx=\int_{\Omega}\operatorname{Pr}_{C}(\frac{p_{0,h}}{\gamma})\chi_{O}\cdot\varphi dx \quad\text{for all}\quad \varphi\in V_{0h}
    \end{aligned}
    \right.
    \end{equation}
    
    and for $n=1,2,\ldots,N-1$, solve 
    \begin{equation}\label{opt_condition4}
    \left\{
    \begin{aligned}
    &z_{n,h}\in V_{0,h}\\
    &\int_{\Omega}\frac{z_{n,h}-z_{n-1,h}}{\Delta t}\varphi dx+\nu\int_{\Omega}\nabla z_{n,h}\cdot\nabla\varphi dx+\int_{\Omega}a_0z_{n,h}\varphi dx=\int_{\Omega}\operatorname{Pr}_{C}(\frac{p_{n,h}}{\gamma})\chi_{O}\cdot\varphi dx \quad\text{for all}\quad \varphi\in V_{0h}
    \end{aligned}
    \right.
    \end{equation}
  
   The strategy for gradient computing advocated here belongs to discretize-then-optimize. Precisely, we first discretize problem and compute the gradient in the discretized setting. Thus the discrete state equation (\ref{opt_condition2}) and discrete adjoint equation (\ref{opt_condition3})-(\ref{opt_condition4}) are strictly in duality which guarantees that the direction $-DJ_h^{\Delta t}(q_h^{\Delta})$ is a discent direction for functional $J_h^{\Delta t}$ at $q_h^{\Delta}$.
  
   \textbf{Remark.} An alternative can be advocated: firstly, derive the adjoint equation to compute the first-order differential of the cost functional in the continuous setting, then discretize the state and adjoint equations simultaneously by some certain numerical scheme, finally compute a discretization of the differential of the cost functional by discretized state solution and adjoint solution. This numerical scheme for gradient computation belongs to optimize-then-discretize. The main problem of this scheme is that the strict duality between the discrete state equation and the discrete adjoint equation may not be preserved. Thus the gradient derived by this scheme may not be the gradient of discretized problem. As a result, the resulting algorithm may not be a descent algorithm and divergence may even appear as discussed in \cite{MR1646758}. 
   
   Finally we conclude this section by giving an implementable algorithm for fully discrete problem (\ref{full_discrete_prob}) that can be regarded as discrete analogue of \textbf{(a)}-\textbf{(d)}.
      
      \begin{algorithm}[htb]
      	\caption{Dual+FRCG} 
      	\hspace*{0.02in} {\bf Step 1:} Give an initial guess $q^0=\{q^0_{h,n}\}_{n=0}^{N-1}$, stopping tolerence constant $tol$ and line search constant $c$. Set iteration number $k=0$.\\
      	\hspace*{0.02in} {\bf Step 2:}
      	Obtain the gradient of objective functional at $q^0$ denoted as $g^0=DJ_h^{\Delta t}(q^0)$ by solving two successive parabolic equations corresponding to $q^0$ (\ref{opt_condition2}) and (\ref{opt_condition3})-(\ref{opt_condition4}).
      	
      	 If $\|g^0\|_{\Delta t}<tol$, then set $\bar{q}=q^0$ as solution and go to \textbf{Step 5}.; otherwise set $d^0=-g^0$. \\
      	\hspace*{0.02in} {\bf Step 3:} Choose the stepsize $\rho_k$ satisfying the following condition:
      	\begin{equation*}
      	J_h^{\Delta t}(q^k+\rho_kd^k)\leq J_h^{\Delta t}(q^k)+ c\cdot\langle DJ_h^{\Delta t},d^k\rangle_{\Delta t}\cdot\rho_k 
      	\end{equation*}\\
      	\hspace*{0.02in} {\bf Step 4:} 
      	Update $q^{k+1}$ and $g^{k+1}$ respectively by
      	\begin{equation*}
      	\begin{aligned}
      	q^{k+1}=q^k+\rho_kd^k\\
      	g^{k+1}=DJ_h^{\Delta t}(q^{k+1})
      	\end{aligned}
      	\end{equation*}
      	
      	If $\|g^{k+1}\|_{\Delta t}<tol$, take $\bar{q}=q^{k+1}$ as solution and go to \textbf{Step 5}.
      	 
      	Otherwise compute 
      	\begin{equation*}
      	\beta_k=\frac{\|g^{k+1}\|^2_{\Delta t}}{\|g^k\|^2_{\Delta t}}
      	\end{equation*}
      	and then update $d^{k+1}$ by
      	\begin{equation*}
      	d^{k+1}=-g^{k+1}+\beta_kd^k
      	\end{equation*}
      	
      	Set $k=k+1$ and return to \textbf{Step 3}.\\
      	\hspace*{0.02in} {\bf Step 5:} Obtain primal problem's solution by formular (\ref{rel_pri_dual}).
      \end{algorithm} 
 
    \section{Second-Order Algorithm Design}\label{SSN}
 
 The regularization parameter $\gamma$ in our problem (\ref{eq3})-(\ref{eq4}) can be set very small. It is evidently that the objective functional in primal problem or dual problem will become increasingly ill-conditioned as $\gamma$ decreases. As a result, for smaller constant $\gamma$, choosing first-order algorithm to solve problem is not suitable. Thus in this section we design a second-order algorithm to solve problem that is based on the dual problem and semismooth Newton computational framework.
 
 \subsection{The discrete optimality system (\ref{opt_sys1})-(\ref{opt_sys2})}
 
 Here we apply the optimize-then-discretize strategy to solve the dual problem (\ref{dual_prob})-(\ref{dual_prob_con}) which means we directly solve the optimality conditions (\ref{opt_sys1})-(\ref{opt_sys2}). Applying the backward-Euler finite difference method for time discretization with $N$ time steps of size $\Delta t=\frac{T}{N}$ and piecewise linear element method for space discretization, gives the following discretized optimality system. 
 \begin{align}
 \left[
 \begin{array}{c:c}
 \mathcal{M} &\mathcal{K}^{\top}\\
  \hdashline
 \mathcal{K} &-\mathcal{M}_1\operatorname{Pr}^{\mathcal{M}}_{\mathcal{C}}(\frac{\cdot}{\gamma})\\
 \end{array}
 \right]
 \left[
 \begin{array}{c}
 \bs z\\
 \bs p
 \end{array}
 \right]&=
 \left[
 \begin{array}{c}
 \mathcal{M}\bs{y_d}\\
 \hdashline
 0
 \end{array}
 \right] \label{opt_s1}\\
  z^0&=0\\
 \widehat{K}p^0-\frac{M}{\Delta t}p^1&=My_d^0\\
 -\frac{M}{\Delta t}z^{N-1}+\widehat{K}z^N&=0\\
  p^N&=0 \label{opt_s2}
 \end{align}
 where $\bs z$, $\bs p$ and $\bs{y_d}$ denote vector corresponding to the state, adjoint and desired state at time-steps $1,2,\ldots,N-1$, and 
 
\begin{equation}
\begin{aligned}
&\mathcal{M}=
\left[\begin{array}{ccccc}
M &  & & & \\
&M &  & & \\
& & \ddots &  & \\
& & & & M 
\end{array}\right]\quad
\mathcal{K}=
\left[\begin{array}{ccccc}
\widehat{K} &  & & & \\
-\frac{M}{\Delta t} &\widehat{K} &  & & \\
& \ddots & \ddots & & \\
& &-\frac{M}{\Delta t} &\widehat{K} 
\end{array}\right]\\
&\widehat{K}=\frac{M}{\Delta t}+\nu K+a_0M\\
&\mathcal{M}_1=
\left[\begin{array}{ccccc}
M_1 &  & & & \\
&M_1 &  & & \\
& & \ddots &  & \\
& & & & M_1 
\end{array}\right]\quad
\mathcal{M}_1\operatorname{Pr}^{\mathcal{M}}_{\mathcal{C}}(\frac{\cdot}{\gamma})
=\left[\begin{array}{ccccc}
M_1\operatorname{Pr}^M_{\mathcal{C}}(\frac{\cdot}{\gamma}) &  & & & \\
&M_1\operatorname{Pr}^M_{\mathcal{C}}(\frac{\cdot}{\gamma}) &  & & \\
& & \ddots &  & \\
& & & &M_1\operatorname{Pr}^M_{\mathcal{C}}(\frac{\cdot}{\gamma}) 
\end{array}\right]
\end{aligned}
\end{equation}

Here, $M$ denotes a finite element mass matrix over the space domain $\Omega$; similarly $M_1$ denotes the finite mass matrix for the domain $\mathcal{O}$ and $K$ a stiffness matrix over $\Omega$. These are defined by
 \begin{equation*}
 \begin{aligned}
 &M=(m_{ij})_{n\times n},\quad m_{ij}=\int_{\Omega} \phi_i\phi_j dx,\\
 &M_1=(m^1_{ij})_{n\times n},\quad m^1_{ij}=\int_{\mathcal{O}} \phi_i\phi_j dx,\\
 &K=(k_{ij})_{n\times n},\quad k_{ij}=\int_{\Omega} \nabla\phi_i\cdot\nabla\phi_j dx.
 \end{aligned}
 \end{equation*}
 
 And $\operatorname{Pr}^M_{\mathcal{C}}(\cdot)$ denotes the projection onto $\mathcal{C}$ w.r.t. the norm $\|\cdot\|_M$. Consequently, $v_h=\operatorname{Pr}^M_{\mathcal{C}}(\omega_h)$ if and only if
 \begin{equation}\label{proj}
 (v_h-\omega_h)^{\top}M(u_h-v_h)\geq 0,\quad\forall u_h\in V_h\cap\mathcal{C}
 \end{equation}
 
 Note that the projection formula (\ref{proj}) cannot be evaluated in a specific manner. To address this problem we consider mass lumping technique, precisely consider $M$ and $M_1$ to be a lumped mass matrix, that is,
 \begin{equation*}
 M=\operatorname{diag}(m_{ii}),\quad m_{ii}=\sum\limits_{j=1}^n\left|\int_{\Omega}\phi_i\phi_jdx\right|
 \end{equation*}
 The $M_1$ can be obtained by replacing $\Omega$ by $\mathcal{O}$. Taking into account this fact that $M$ is a diagonal matrix, the formular (\ref{proj}) can be evaluated specifically, precisely, there holds
 \begin{equation*}
 \operatorname{Pr}^M_{\mathcal{C}}(\cdot)=\operatorname{Pr}_{\mathcal{C}}(\cdot)
 \end{equation*}
 
 \textbf{Remark.} From now on, we replace mass matrices $M$ and $M_1$ by lumped mass matrices. It is reasonable both from computational and theoretical points of view. Mass lumping is a standard tool for the numerical solution of time-dependent pde \cite{MR2249024}. Besides, many algorithms for solving pde optimal control problems are designed based on mass lumping technique \cite{MR3023467,MR3504560,MR2114385}. Furthermore, some rigorous results corresponding to the error analysis about mass lumping technique applied to optimal control problems, we refer to \cite{MR3663005} for more detail.
 
 \subsection{The active-set Newton method}
 In the following, we derive an active-set Newton type method for the solution of discrete optimality system (\ref{opt_s1})-(\ref{opt_s2}). 
 
 Since the system is of some special structure, we observe that as long as the nonlinear equations (\ref{opt_s1}) is solved then the whole optimality system is solved. Thus the main difficulty is converted into how to solve (\ref{opt_s1}) efficiently. 
 
 Let us denote (\ref{opt_s1}) as
 \begin{equation}\label{opte_s1}
 F(\bs z,\bs p)=\left[
 \begin{array}{c}
 \mathcal{M}(\bs z-\bs{y_d})+\mathcal{K}^{\top}\bs p\\
 \mathcal{K}\bs z -\mathcal{M}_1\operatorname{Pr}_{\mathcal{C}}(\frac{\bs p}{\gamma})\\
 \end{array}
 \right]=0
 \end{equation}
 
 \begin{lemma}
 	The equations (\ref{opt_s1}) has a unique solution. 
 \end{lemma}
 
 {\emph{Proof.}} 
 The equations (\ref{opt_s1}) can be equivalently represented as:
 \begin{align}
 &\mathcal{M}\bs z+\mathcal{K}^{\top}\bs p=\mathcal{M}\bs{y_d} \label{non_1}\\
 &-\mathcal{M}_1\operatorname{Pr}_{\mathcal{C}}(\frac{\bs p}{\gamma})-\mathcal{K}\mathcal{M}^{-1}\mathcal{K}^{\top}\bs p=-\mathcal{K}\bs {y_d}\label{non_2}
 \end{align}
 
 Since the matrix $\mathcal{M}_1$ is diagonal and $\mathcal{K}$ is a full rank matrix, thus the operator $\mathcal{M}_1\operatorname{Pr}_{\mathcal{C}}(\frac{\cdot}{\gamma})+\mathcal{K}\mathcal{M}^{-1}\mathcal{K}^{\top}$ is maximal monotone. This means there must exist one and only one $\bs p^*$ to satisfy equation (\ref{non_2}). 
 
 The equation (\ref{non_1}) can be tranformed into the following equation because $\mathcal{M}$ is inversable.
 \begin{equation*}
 \bs z=\bs{y_d}-\mathcal{M}^{-1}\mathcal{K}^{\top}\bs p
 \end{equation*}
 
 Thus there must exist a unique solution $(\bs z^*,\bs p^*)$ such that the nonlinear equations (\ref{opt_s1}) hold.
 $\hfill\qedsymbol$ 
 
 The nonlinearity and nonsmoothness of the function $F$ defined in (\ref{opte_s1}) are gathered in the second diagonal block containing the projection operators. This fact suggests that we can use the generalized Jacobian to construct a "semismooth" Newton scheme.
 
 Given the $k$-th iteration point $(\bs z^k,\bs p^k)$, precisely 
 $\bs z^k=[z^k_1;z^k_2;\ldots;z^k_{N-1}]$ and $\bs p^k=[p^k_1;p^k_2;\ldots;p^k_{N-1}]$
 
 Let $\mathcal{A}^k_i$ denotes the current active set corresponding to $p^k_i$,
 \begin{equation*}
 \mathcal{A}^k_i=\left\{j\mid\frac{( p^k_i)_j}{\gamma}\in\mathcal{C}\right\}
 \end{equation*}
 and let $\Pi_{\mathcal{A}^k_i}$ denotes a diagonal binary matrix with nonzero entries in $\mathcal{A}^k_i$. We define $\Pi_k$ by following formular:
 \begin{equation*}
 \Pi_k=
 \left[\begin{array}{ccccc}
 	\Pi_{\mathcal{A}^k_1} &  & & & \\
 	&\Pi_{\mathcal{A}^k_2} &  & & \\
 	& & \ddots &  & \\
 	& & & &\Pi_{\mathcal{A}^k_{N-1}}  
 \end{array}\right]
 \end{equation*}
  
 Then the generalized Jacobian matrix of $F$ (\ref{opte_s1}) at $(\bs z^k,\bs p^k)$ can be given by
 \begin{equation}\label{Gen_J}
 F^{\prime}(\bs z^k,\bs p^k)=
 \left[\begin{array}{cc}
 \mathcal{M} &\mathcal{K}^{\top}\\
 \mathcal{K} &-\frac{\mathcal{M}_1\Pi_k}{\gamma}
 \end{array}
 \right]
 \end{equation}
 Using the generalized Jacobian matrix above, the following formular conceptually is the semismooth Newton iteration applied to original nonlinear system (\ref{opte_s1}):
 \begin{equation}\label{Newton_sys}
 F(\bs z^k,\bs p^k)+F^{\prime}(\bs z^k,\bs p^k)\left[
 \begin{array}{c}
 \bs z^{k+1}-\bs z^k\\
 \bs p^{k+1}-\bs p^k
 \end{array}
 \right]=\left[\begin{array}{c}
  \mathcal{M}\bs{y_d}\\
 0
 \end{array}
 \right]
 \end{equation}
 
 Since $F^{\prime}(\bs z^k,\bs p^k)$ must be inversable, the Newton equation (\ref{Newton_sys}) exist unique solution.
 
 Finally we conclude this subsection by the following numerical scheme for solving (\ref{opte_s1})
 
 \begin{algorithm}[htb]
 	\caption{Dual+SSN}
 	\hspace*{0.02in} {\bf Step 1:}  Give an initial guess $(\bs z^0,\bs p^0)$, stopping tolerence constant $tol$. Set iteration number $k=0$. \\
 	\hspace*{0.02in} {\bf Step 2:}
 	Construct Newton equation (\ref{Newton_sys}) by $k$-th iteration point $(\bs z^k,\bs p^k)$.\\
 	\hspace*{0.02in} {\bf Step 3:}
 	Update $(\bs z^{k+1},\bs p^{k+1})$ by	
 	solving Newton equation obtained at (\textbf{Step 2}).\\
 	\hspace*{0.02in} {\bf Step 4:} If the following inequality holds
 	\begin{equation*}
 	\left\|\mathcal{K}\bs z^{k+1} -\mathcal{M}_1\operatorname{Pr}_{\mathcal{C}}(\frac{\bs p^{k+1}}{\gamma})\right\|\leq tol
 	\end{equation*}
 	then take $(\bs z^{k+1},\bs p^{k+1})$ as solution and go to \textbf{Step 5}. 
 	
 	Otherwise set $k=k+1$ and return to \textbf{Step 2}.\\
 	\hspace*{0.02in} {\bf Step 5:} Obtain primal problem's solution by formular (\ref{rel_pri_dual}).
 \end{algorithm}

 \textbf{Remark.} The above numerical scheme derived by us has been proved to be of locally superlinear convergence rate and locally convergence \cite{MR2839219}. For the locally convergence result, there exist some globalization methods and these can be directly embedded in our algorithmic design. Besides, for each Newton system there also exist some results suggesting that it can be solved inexactly, furthermore the locally convergence rate is still retained. These results are important but beyond the scope of our discussion, we refer to \cite{MR1354652} for more detail. 
 
 \subsection{Solving the Newton equation}
 
 In the following, we consider how to solve each Newton equation efficiently. It is obviously that each linear equation is ill-conditioned and of large-scale thus this discussion is a must.  
 
 For simplicity, we introduce the following notation to represent system (\ref{Newton_sys}):
 \begin{equation*}
 \begin{aligned}
 &\Delta \bs z^k:=\bs z^{k+1}-\bs z^k,\quad \Delta \bs p^k:=\bs p^{k+1}-\bs p^k\\
 &\bs d^k:=\left[\begin{array}{c}
 \mathcal{M}\bs{y_d}\\
 0
 \end{array}
 \right]-F(\bs z^k,\bs p^k)
 \end{aligned}
 \end{equation*}
 Thus the original system (\ref{Newton_sys}) can be represented equivalently as
 \begin{equation}\label{line_1}
 F^{\prime}(\bs z^k,\bs p^k)\left[
 \begin{array}{c}
 \Delta \bs z^k\\
 \Delta \bs p^k
 \end{array}
 \right]=\bs d^k
 \end{equation}
 
 For the matrix $F^{\prime}(\bs z^k,\bs p^k)$, it can be factorized as
 \begin{equation}
 F^{\prime}(\bs z^k,\bs p^k)=\underbrace{\left[
 \begin{array}{cc}
 I & 0\\
 \mathcal{K}\mathcal{M}^{-1} & I
 \end{array}
 \right]}_{L}\cdot
 \underbrace{\left[
 \begin{array}{cc}
 \mathcal{M} & 0\\
 0 & -(\frac{\mathcal{M}_1\Pi_k}{\gamma}+\mathcal{K}\mathcal{M}^{-1}\mathcal{K}^{\top})
 \end{array}
 \right]}_{\operatorname{blkdiag}(\mathcal{M},-C_k)}\cdot
 \underbrace{\left[
 \begin{array}{cc}
 I & \mathcal{M}^{-1}\mathcal{K}^{\top}\\
 0 & I
 \end{array}
 \right]}_{L^{\top}}
 \end{equation}

 Hence, the procedure of solving linear equation (\ref{line_1}) can be summarized as following:
  \begin{enumerate}
 	\item[\textbf{(a)}] Solve linear system: $L\widehat{\bs d}^k=\bs d^k$
 	\item[\textbf{(b)}] Solve linear system: $\operatorname{blkdiag}(\mathcal{M},-C_k)\bar{\bs d}^k=\widehat{\bs d}^k$
 	\item[\textbf{(c)}] Solve linear system: $L^{\top}\left[
 	\begin{array}{c}
 	\Delta \bs z^k\\
 	\Delta \bs p^k
 	\end{array}
 	\right]=\bar{\bs d}^k$
 	
 \end{enumerate}
 
 Since $M$ is a diagonal matrix, step \textbf{(a)} and step \textbf{(c)} are easy to compute. For step \textbf{(b)}, it is much more difficult to approach mainly because stiffness matrix $K$ is ill-conditioned and it appear in $C_k$. It motivates us to design a preconditioner for solving linear equation obtained by step \textbf{(b)}. Our main idea is to approximate the second block of matrix $C_k$, precisely the Schur complement of $F^{\prime}(\bs z^k,\bs p^k)$. This procedure is mainly inspired by \cite{MR3023467,MR2979813}.
 
 We define the following factorized approximation of $C_k$:
 \begin{equation*}
 \mathbb{C}_k=(\mathcal{K}+\frac{\mathcal{M}^{\frac{1}{2}}\mathcal{M}_1^{\frac{1}{2}}}{\sqrt{\gamma}}\Pi_k)\mathcal{M}^{-1}(\mathcal{K}+\frac{\mathcal{M}^{\frac{1}{2}}\mathcal{M}^{\frac{1}{2}}_1}{\sqrt{\gamma}}\Pi_k)^{\top}
 \end{equation*}
 
 \textbf{Remark.} Our approximation mainly uses the fact that $\mathcal{M}$, $\mathcal{M}_1$ and $\Pi_k$ are diagonal, meanwhile the element of $\Pi_k$ is binary. 
 
 Then we should analyze the quality of the proposed preconditioner $\mathbb{C}_k$, precisely, deriving the spectral property of $\mathbb{C}_k^{-1}C_k$. In order to do this, we need to prove the following Lemma firstly.
 
 \begin{lemma}\label{lemma_pd}
 	The matrix $\mathcal{K}+\mathcal{K}^{\top}$ is positive definite.
 \end{lemma}
 {\emph{Proof.}}
 For simplicity, we define matrix $E_1$ as follows:
 \begin{equation*}
 E_1=\left[\begin{array}{ccccc}
 0 &  & & & \\
 1&0 &  & & \\
 & 1& 0 &  & \\
& & \ddots &\ddots & \\
 & & & 1&0 
 \end{array} 
 \right]
 \end{equation*}
 
 We observe that matrix $\mathcal{K}$ can be represented as
 \begin{equation*}
 \begin{aligned}
 \mathcal{K}&=I\otimes\widehat{K}-E_1\otimes \frac{M}{\Delta t}\\
 &=(I-E_1)\otimes\frac{M}{\Delta t}+I\otimes(\nu K+a_0M)
 \end{aligned}
 \end{equation*}
 
 Thus there holds:
 \begin{equation*}
 \mathcal{K}+\mathcal{K}^{\top}=(2I-E_1-E_1^{\top})\otimes\frac{M}{\Delta t}+I\otimes(2\nu K+2a_0M)
 \end{equation*}
 
 The matrix $2I-E_1-E_1^{\top}$ is positive definite. Because of the nonnegative coefficient $\nu$ and $a_0$, the matrix $2\nu K+2a_0M$ is also positive semidefinite. The conclusion holds due to the property of kronecker product.
 $\hfill\qedsymbol$
 
 Then we give the conclusion of spectral property.
 \begin{theorem}\label{th_spectral}
 	Let $\lambda$ be an eigenvalue of $\mathbb{C}_k^{-1}C_k$. Then there holds:
 	\begin{equation*}
 	\frac{1}{2}\leq\lambda\leq\zeta^2+(1+\zeta)^2
 	\end{equation*}
 	with $\zeta=\sqrt{\gamma}\|(\sqrt{\gamma}I+\mathcal{M}^{\frac{1}{2}}\mathcal{K}^{-1}\mathcal{M}^{\frac{1}{2}}\Pi_k)^{-1}\|$
 	
 	Furthermore, there holds for $\gamma\to 0^+$, $\lambda$ can be bounded by a constant independent of $\gamma$.  
 \end{theorem}
 {\emph{Proof}.} See the Appendix \ref{Appendix} for the proof. $\hfill\qedsymbol$
 
 Then we can specify how to solve linear equation $C_k\bar{\bs d}^k=\widehat{\bs d}^k$ (here for simplicity we still denote as $\bar{\bs d}^k$, $\widehat{\bs d}^k$ but these are different  from the notations occur in Step \textbf{(b)}).  
 
 It is obviously that the matrix $C_k$ is positive definite, thus we can use PCG to solve this large-scale linear equation. As $C_k$ is ill-conditioned, the preconditioner for it is chosen as $\mathbb{C}_k$. Because in each PCG iteration, a large-scale linear equation w.r.t. preconditioner $\mathbb{C}_k$ need to be solved, we discuss how to solve it efficiently as follows. 
 
 Notice that the matrix $(\mathcal{K}+\frac{\mathcal{M}^{\frac{1}{2}}\mathcal{M}_1^{\frac{1}{2}}}{\sqrt{\gamma}}\Pi_k)$ is a block lower triangular matrix, thus we can solve the equation w.r.t. $\mathbb{C}_k$ by a forward sweep w.r.t. $(\mathcal{K}+\frac{\mathcal{M}^{\frac{1}{2}}\mathcal{M}_1^{\frac{1}{2}}}{\sqrt{\gamma}}\Pi_k)$ and a backward sweep w.r.t. $(\mathcal{K}+\frac{\mathcal{M}^{\frac{1}{2}}\mathcal{M}_1^{\frac{1}{2}}}{\sqrt{\gamma}}\Pi_k)^{\top}$.  
 
 In both forward sweep procedure and backward sweep procedure, we need to solve linear equation w.r.t. $\widehat{K}+\frac{M^{\frac{1}{2}}M_1^{\frac{1}{2}}}{\sqrt{\gamma}}\Pi_{\mathcal{A}^k_i}$. We note that it is still ill-conditioned and directly applying its inverse may be not feasible. Hence, for a practical algorithm, we choose a multigrid V-cycles associated with it to approximate its inverse matrix.
 
 \textbf{Remark.} We note that each iteration the Newton equation preconditioner designed by us is of the same structure, and the dimension of Newton equation is also invariant. Besides there is only one variant block in Newton equation, which is different from \cite{MR3504560}. Because of these characteristics in each iteration, our designed algorithm is implementable.

\section{Numerical Result}	
	In this section we numerically verify the efficiency of our designed first-order and second-order algorithms. We compare our algorithm called 'Dual+FRCG' and 'Dual+SSN' with inexact ADMM called 'In-ADMM'. 
	 Our codes were written in MATLAB R2019b and all numerical experiments were conducted on a computer with the process, Inter(R) Core(TM) i7-7660U CPU at 2.50GHz, and with a 32.00-GB RAM.
	 
	 To test our proposed algorithm: 'Dual+FRCG' we set the stopping criterion as
	 \begin{equation*}
	 \frac{\Delta t\sum\limits_{n=0}^{N-1}\|g^k_n\|^2}{\Delta t\sum\limits_{n=0}^{N-1}\|g^0_n\|^2}\leq tol^2
	 \end{equation*}
	 
	 For 'Dual+SSN' we set the stopping criterion as
	 \begin{equation}\label{SSN_tol}
	 \left\|\mathcal{K}\bs z^{k} -\mathcal{M}_1\operatorname{Pr}_{\mathcal{C}}(\frac{\bs p^{k}}{\gamma})\right\|\leq tol
	 \end{equation}
	 
	 For 'Dual+FRCG', we choose $tol=10^{-4}$ and initial values are chosen as $\bs q=0$. The parameter for inexact line-search condition (\ref{inexact_line}) is set by $c=0.4$ 
	 
	 For 'Dual+SSN', we choose $tol=10^{-4}$ and initial values are chosen as $\bs z=0$ and $\bs p=0$. The large-scale linear equation obtained each step is solved by 'pcg' solver in MATLAB and tolerance is set by $10^{-6}$. Futhermore solving each linear equation w.r.t. preconditioner as descibed in previous section, we should do forward sweep and backward sweep. We choose multigrid V-cycles to solve linear systems appearing in both sweep procedure, and this implementation is based on the iFEM package developed in \cite{chen2008ifem}.
	 
	 For 'In-ADMM', the primal residual and dual residual are denoted as $\pi_s$ and $d_s$ respectively. The stopping criteria for 'In-ADMM' for all numerical experiments is set by
	 \begin{equation*}
	 \max\{\pi_s,d_s\}\leq tol
	 \end{equation*}
	 The constant $tol=10^{-4}$ and initial values are set as $u=0$, $z=0$ and $\lambda=0$. For the constant $\sigma$ defined in inexactness criterion, we choose $\sigma=0.99\frac{\sqrt{2}}{\sqrt{2}+\sqrt{\beta}}$. For more details about 'In-ADMM', we refer reader to \cite{song2020implementation}.
	 
	 Besides for each linear system arising at each time step of the discretized parabolic equations in 'In-ADMM' or 'Dual+FRCG', they are also solved by multigrid V-cycles.
	 
	 In addition, we define the relative distance "RelDis" and the objective functional value "Obj" as:
	 \begin{equation*}
	 \text{ReDis}=\frac{\|y-y_d\|^2_{L^2(Q)}}{\|y_d\|^2_{L^2(Q)}},\quad \text{Obj}=\frac{1}{2}\|y-y_d\|^2_{L^2(Q)}+\frac{\gamma}{2}\|u\|^2_{L^2(\mathcal{O})}
	 \end{equation*}
	  
	  	For all our numerical experiments the space mesh size $h$ and time steps $\Delta t$ are set as $h=\Delta t=2^{-i}$ with $i=4,5,6,7,8$. In all numerical table, notation 'Iter' denotes the total out-layer iteration number while 'Mean/Max CG' denotes the average and maximum steps of the inner CG method of inexact ADMM, for simplicity, 'Mean/Max CG' also denotes the inner PCG method of semismooth Newton method. One should note that each iteration implemented by ADMM type method is 2 layer-nested, while 'Dual+FRCG' is one layer and 'Dual+SSN' is also 2 layer-nested. Because 'Dual+SSN' requires solving large-scale Newton equation each step and 'In-ADMM' needs solve an uncontrained subproblem.
	  
	\textbf{Example 1.} We consider the following example with a known exact solution. The model is adapted from \cite{MR2891918}.
	\begin{equation*}
	\begin{aligned}
	\min\limits_{u \in \mathcal{C}, y \in L^{2}(Q)}\quad & \frac{1}{2} \iint_{Q}\left|y-y_{d}\right|^{2} d x d t+\frac{\gamma}{2} \iint_{Q}|u|^{2} d x d t \\
	\mbox{ s.t. }\quad & \left\{\begin{array}{ll}
	\frac{\partial y}{\partial t}-\Delta y=f+u, & \text { in } \Omega \times(0, T) \\
	y=0, & \text { on } \Gamma \times(0, T) \\
	y(0)=\varphi
	\end{array}\right.
	\end{aligned}
	\end{equation*}
	
	with $\Omega=(0,1)^{2}, T=1$. The function $f \in L^{2}(Q)$ is a source term that helps us construct the exact solution without affection to the numerical implementation. We further let
	\begin{equation*}
	\left\{\begin{array}{l}
	y^*=(1-t) \sin \pi x_{1} \sin \pi x_{2}, p^*=\gamma(1-t) \sin 2 \pi x_{1} \sin 2 \pi x_{2}, u^*=\min \left(b, \max \left(a,-\frac{p^*}{\gamma}\right)\right) \\
	f=-u^*+\frac{\partial y^*}{\partial t}-\Delta y^*, \quad y_{d}=y^*+\frac{\partial p^*}{\partial t}+\Delta p^*, \quad \varphi=\sin \pi x_{1} \sin \pi x_{2}
	\end{array}\right.
	\end{equation*}
	
	Then, it is obviously that $(u^*,y^*)$ is the optimal solution of the problem. The control admissible set is set as
	\begin{equation*}
	\mathcal{C}=\{v|v\in L^2(Q), -0.5\leq v(t,x)\leq 0.5\text { a.e. in } Q \}
	\end{equation*}
	
	For ADMM type method, it is well-known that the choice parameter $\beta$ is important for the numerical behavior. Here we choose $\beta=3$ as disussed in \cite{song2020implementation}. 
	
\begin{figure}[H]
		\centering
		\begin{minipage}[t]{0.5\textwidth}
			\centering
			\includegraphics[width=6.5cm]{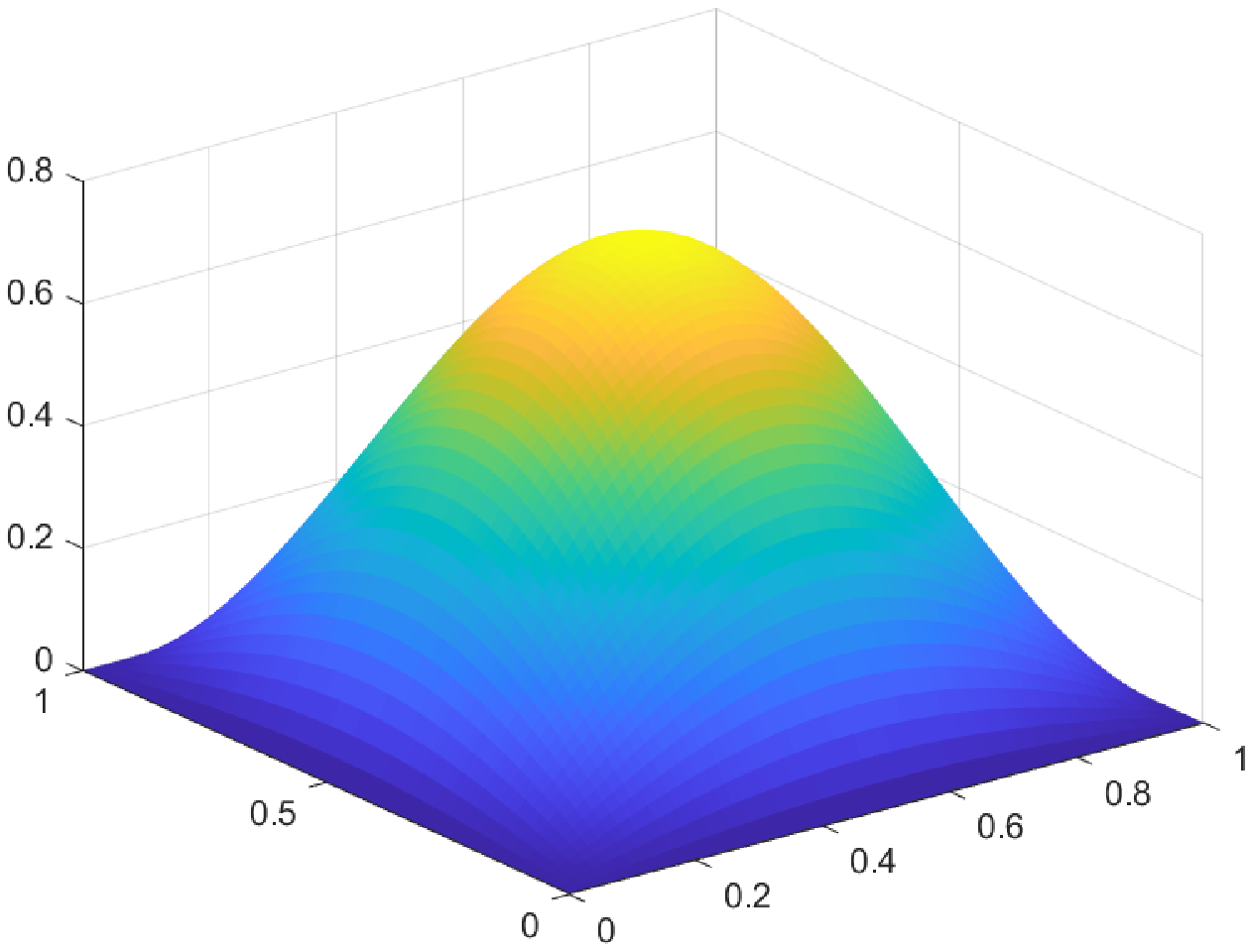}	
		\end{minipage}
		\begin{minipage}[t]{0.49\textwidth}
			\centering
			\includegraphics[width=6.5cm]{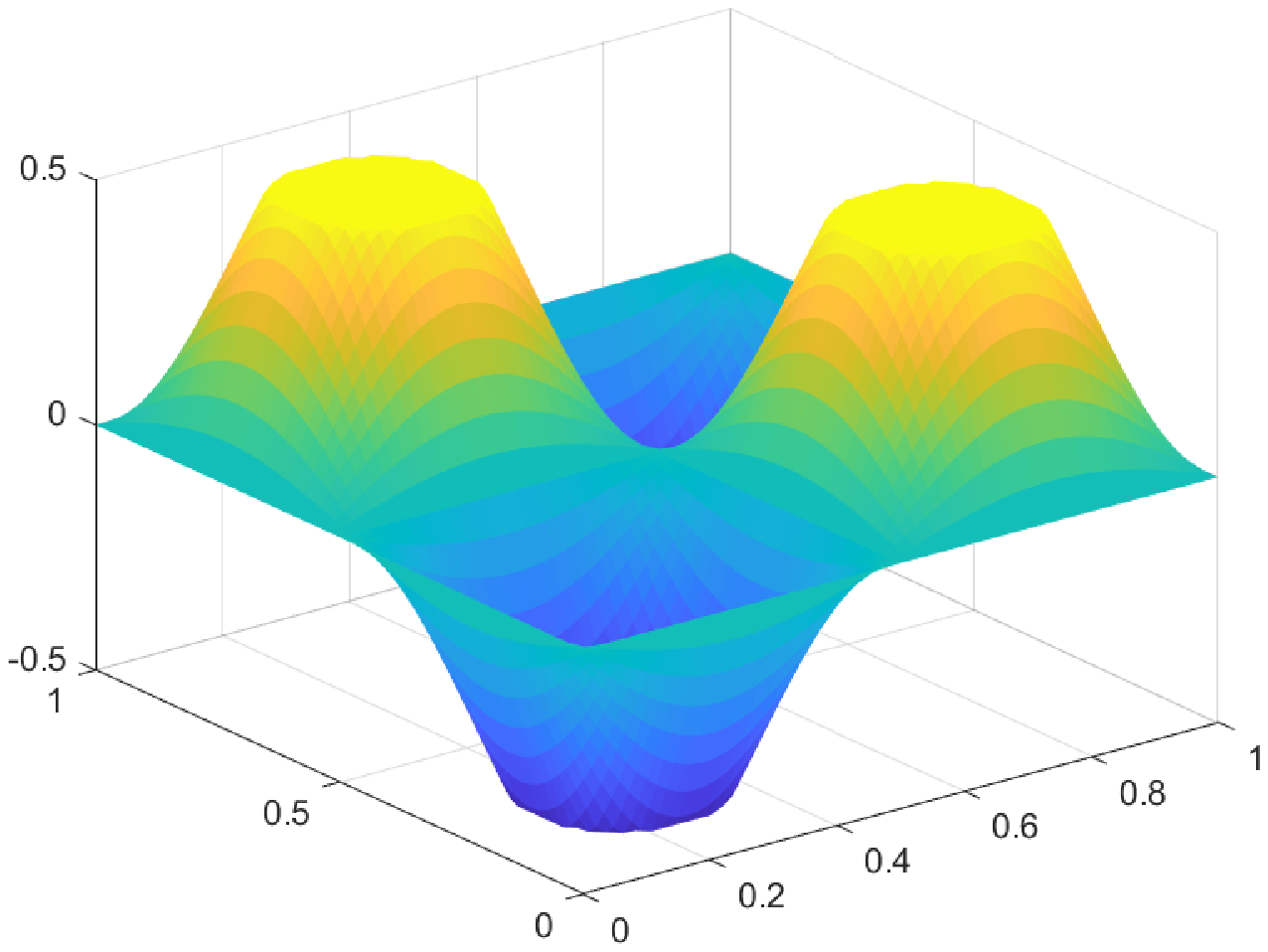}
		\end{minipage}
		\caption{Numercial solution (Dual+FRCG) y (left) and u (right) at t=0.25 for Example 1}
		\label{dual_sol}
	\end{figure}

\begin{figure}[H]
		\centering
		\begin{minipage}[t]{0.5\textwidth}
			\centering
			\includegraphics[width=6.5cm]{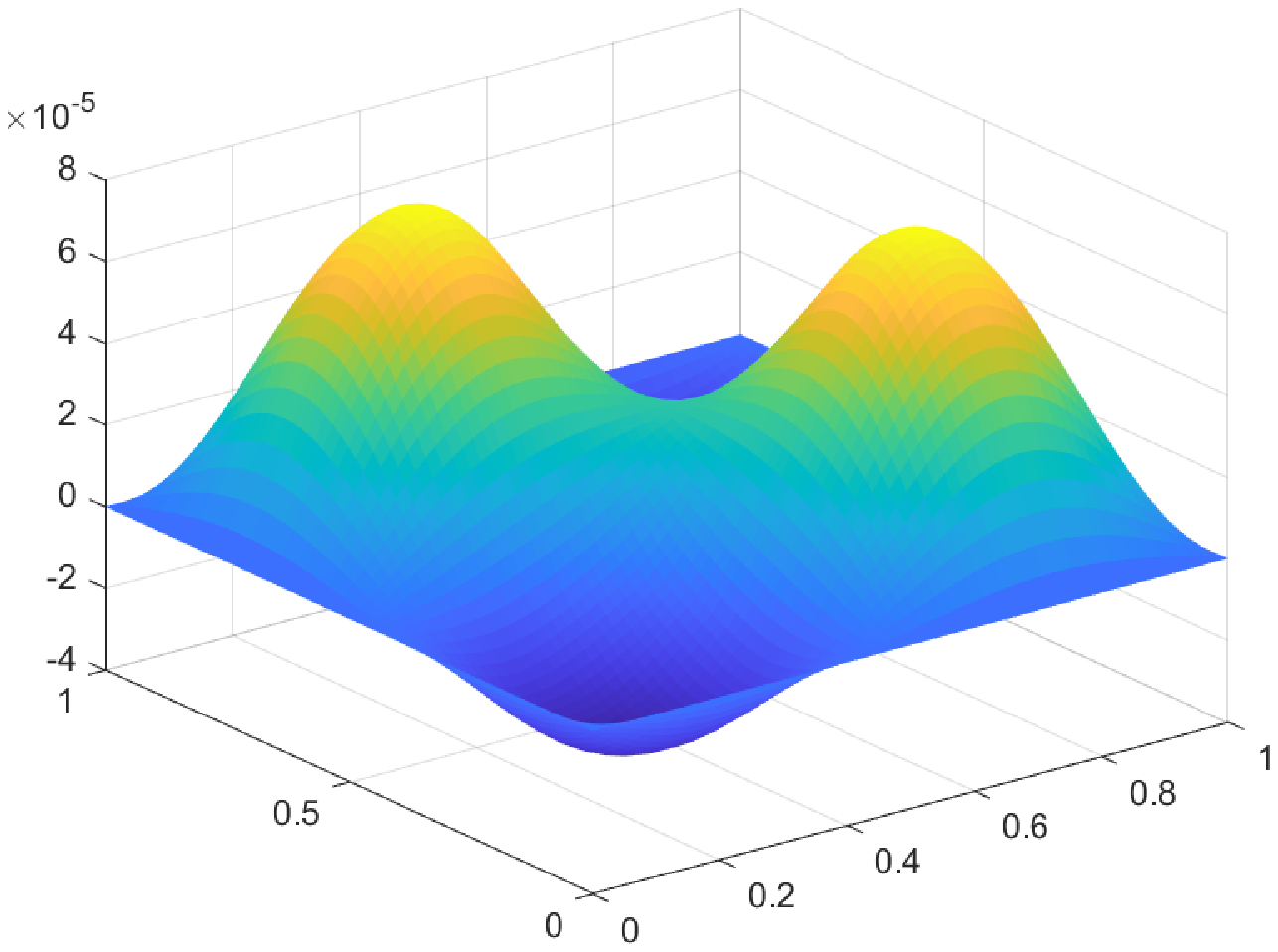}
		\end{minipage}
		\begin{minipage}[t]{0.49\textwidth}
			\centering
			\includegraphics[width=6.5cm]{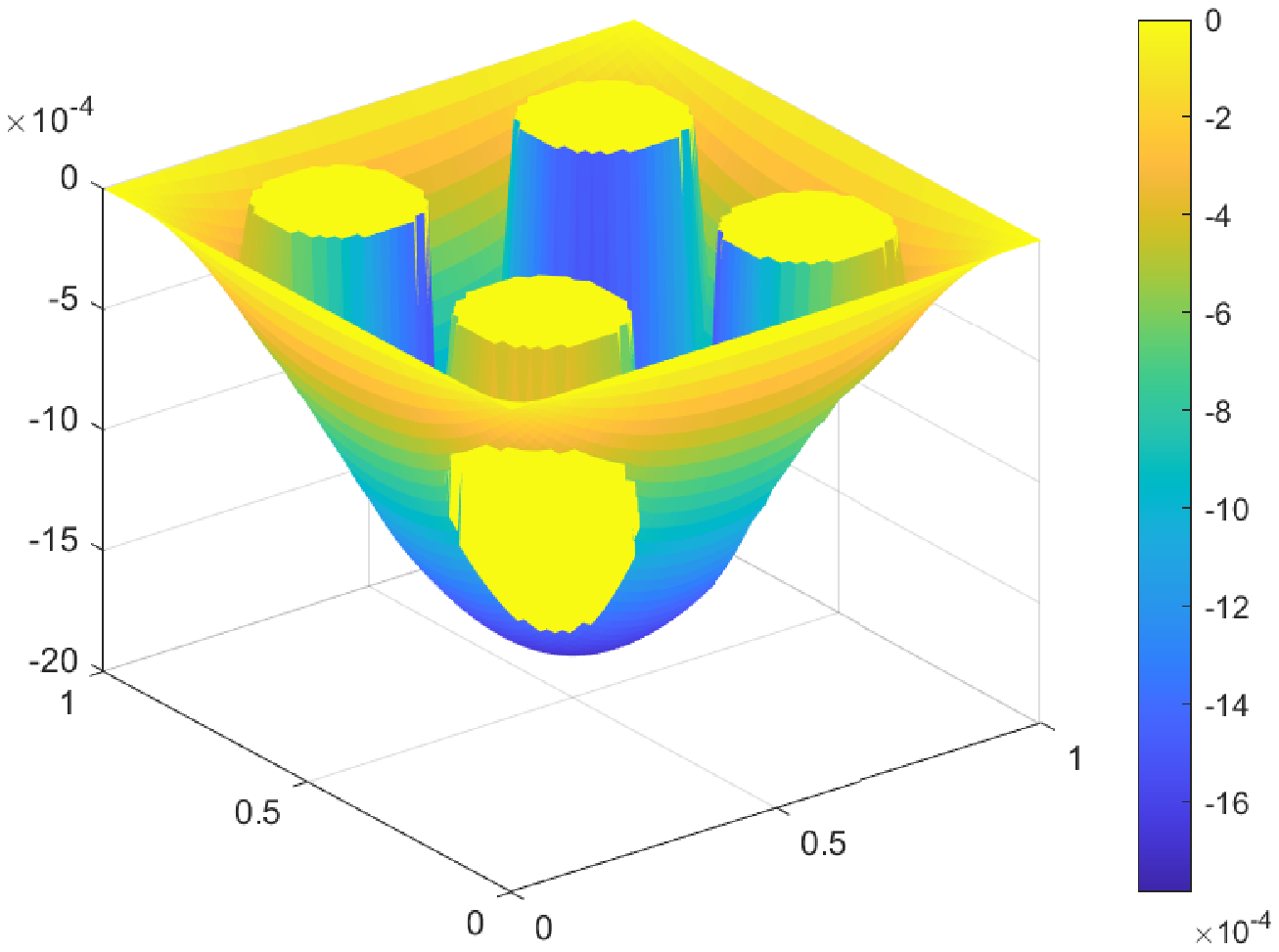}
		\end{minipage}
		\caption{Errors (Dual+FRCG) $y-y^*$ (left) and $u-u^*$ (right) at t=0.25 for Example 1}
		\label{error_dual}
	\end{figure}

	In addition, the state variable $y$ and control variable $u$ obtained by our method (Dual+FRCG) and errors $y-y^*$, $u-u^*$ at $t=0.25$ with $h=\Delta t=2^{-6}$ are depicted in Figure \ref{dual_sol} and \ref{error_dual} respectively.
	
	Firstly, we choose $\gamma=10^{-3}$ and compare our designed first-order algorithm 'Dual+FRCG' with 'In-ADMM'. We note that regularization constant in this case is not close to zero, hence the objective functional is not very 'ill-conditioned'.
	
	\begin{table}[H]
		\caption{Numerical Comparsion of 'In-ADMM' and 'Dual+FRCG' when $\gamma=10^{-3}$}
		\centering
		\begin{tabular}{| l |c|c | c| c|c | c|}
			\hline 
			Mesh &Algorithm &Iter& Mean/Max CG &CPU Time(sec)&Obj &RelDis \\
			\hline
			\multirow{2}{*}{$2^{-4}$} &Dual+FRCG & 9 &\textemdash& 5.41   & $3.28\times 10^{-4}$ & $6.42\times 10^{-3}$ \\
			\cline{2-7}
			\multirow{2}{*}{~} &In-ADMM & 26&1/1 & 7.77  & $3.02\times 10^{-4}$ &$6.43\times10^{-3}$ \\
			\hline
				\multirow{2}{*}{$2^{-5}$} &Dual+FRCG & 8 &\textemdash& 12.9   & $3.15\times 10^{-4}$ & $6.43\times10^{-3}$\\
			\cline{2-7}	
			\multirow{2}{*}{~} &In-ADMM & 26&1/1 & 23.61  & $3.01\times 10^{-4}$ &$6.43\times10^{-3}$ \\
			\hline
			\multirow{2}{*}{$2^{-6}$} &Dual+FRCG & 7 &\textemdash& 62.07   & $3.08\times 10^{-4}$ & $6.43\times10^{-3}$\\
			\cline{2-7}
			\multirow{2}{*}{~} &In-ADMM & 26&1/1 & 165.67  & $3.01\times 10^{-4}$ &$6.43\times10^{-3}$ \\
			\hline
			\multirow{2}{*}{$2^{-7}$} &Dual+FRCG & 9 &\textemdash&896.31   & $3.05\times 10^{-4}$ & $6.43\times10^{-3}$\\
			\cline{2-7}
			\multirow{2}{*}{~} &In-ADMM & 26&1/1 & 1819.86  & $3.01\times 10^{-4}$ &$6.43\times10^{-3}$ \\
			\hline
			\multirow{2}{*}{$2^{-8}$} &Dual+FRCG & 6 &\textemdash&4012.00   & $3.04\times 10^{-4}$ & $6.43\times10^{-3}$\\
			\cline{2-7}
			\multirow{2}{*}{~} &In-ADMM & 26&1/1 & 11009.91  & $3.01\times 10^{-4}$ &$6.43\times10^{-3}$ \\
			\hline
		\end{tabular}
	\label{numerical_table1}
	\end{table}
	
	From Table \ref{numerical_table1}, we observe that at each inexact ADMM iteration, the inner iteration number is just one, that means inexact ADMM behaves very efficient at each iteration. This fact demonstrates that $\beta=3$ is a good choice. But we also notice that 'Dual+FRCG' converge much faster than 'In-ADMM'. The main reason is that our method considers parabolic PDE constraint together with control box constraints while inexact ADMM considers them seperately. Besides total CPU time is also interesting, our algorithm needs to solve two sets of elliptic equations (each set totally include $N=T/\Delta t$) while inexact ADMM involves solving two parabolic PDE in out-layer and each iteration in inner-layer requires solve two parabolic PDE, precisely each iteration inexact ADMM needs solve at least four parabolic PDE. Thus we conclude that when objective functional behaves not very 'ill-conditioned', 'Dual+FRCG' is a good choice, at least for this problem.
	
  Since 'Dual+FRCG'  belongs to first-order algorithm, it is necessary to verify if the solution obtained is close to exact solution. In other words, whether or not it is still the discretization error that dominates the main part of the total error when applying 'Dual+FRCG' solving problem. 
  
  In Table \ref{numerical_error}, we report the $L^2$-error of the solution obtained by 'In-ADMM' and 'Dual+FRCG'. It is clear that when 'Dual+FRCG' is applied to solve the problem, the overall error of $u$ and $y$ are both dominated by the discretization error. This also validates the conclusion in \cite{MR2407017} that the error order of the time discretization is $O(\Delta t)$ and this estimate may
  dominate the magnitude of the total error. Although the error of 'In-ADMM' is smaller than 'Dual+FRCG', we notice that the value of error is very close and our designed first-order algorithm converge much faster. 
  \begin{table}[H]
  	\caption{Numerical errors comparsion of 'In-ADMM' and 'Dual+FRCG' with $\gamma=10^{-3}$}
  	\centering
  	\begin{tabular}{| l |c | c| c|c | c|}
  		\hline 
  		error &Algorithm &$h=\Delta t=2^{-5}$&$h=\Delta t=2^{-6}$&$h=\Delta t=2^{-7}$&$h=\Delta t=2^{-8}$ \\
  		\hline
  		\multirow{2}{*}{$\|u-u^*\|_{L^2(\mathcal{O})}$} &Dual+FRCG &$3.27\times10^{-3} $ &$9.78\times10^{-4}$ &$4.31\times10^{-4}$ & $1.42\times 10^{-4}$  \\
  		\cline{2-6}
  		\multirow{2}{*}{~} &In-ADMM &$3.27\times10^{-3} $ &$8.25\times10^{-4}$ &$2.09\times10^{-4}$ & $8.35\times 10^{-5}$  \\  		
  		\hline
  		\multirow{2}{*}{$\|y-y^*\|_{L^2(Q)}$} &Dual+FRCG &$7.93\times10^{-5}$ &$1.98\times10^{-5}$ & $5.25\times 10^{-6}$ &$1.34\times10^{-6} $  \\
  		\cline{2-6}
  		\multirow{2}{*}{~} &In-ADMM  &$7.80\times10^{-5}$ &$1.95\times10^{-5}$ & $4.97\times 10^{-6}$ &$1.31\times10^{-6} $ \\
  	\hline
  	\end{tabular}
  	\label{numerical_error}
  \end{table}

  Then we set constant $\gamma=10^{-5}$ and compare our designed second-order algorithm 'Dual+SSN' with 'In-ADMM' as following table. For the case $\gamma$ is small, the objective functional of problem is 'ill-conditioned' thus directly using first-order algorithm is not very suitable. And ADMM-type method is still implementable, due to the fact that it considers augmented Lagrangian function at each iteration that is much 'better-conditioned' than original function. But the convergence rate of ADMM-type method may be not good.  

   \begin{table}[H]
   	\caption{Numerical Comparsion of 'In-ADMM' and 'Dual+SSN' when $\gamma=10^{-5}$}
   	\centering
   	\begin{tabular}{| l |c|c | c| c|c | c|}
   		\hline 
   		Mesh &Algorithm &Iter& Mean/Max CG &CPU Time(sec)&Obj &RelDis \\
   		\hline
   		\multirow{2}{*}{$2^{-4}$} &Dual+SSN & 4 &11.5/13& 5.76   & $3.43\times 10^{-7}$ & $6.68\times 10^{-7}$ \\
   		\cline{2-7}
   		\multirow{2}{*}{~} &In-ADMM & 25&5.68/7 & 28.62  & $3.40\times 10^{-7}$ &$6.47\times10^{-7}$ \\
   		\hline
   		\multirow{2}{*}{$2^{-5}$} &Dual+SSN & 4 &12.5/14& 35.16   & $3.41\times 10^{-7}$ & $6.47\times 10^{-7}$ \\
   		\cline{2-7}
   		\multirow{2}{*}{~} &In-ADMM & 22&6.00/7 & 95.67  & $3.41\times 10^{-7}$ &$6.47\times10^{-7}$ \\
   		\hline
   	 \multirow{2}{*}{$2^{-6}$} &Dual+SSN & 4 &13.25/15& 380.03   & $3.41\times 10^{-7}$ & $6.47\times 10^{-7}$ \\
   	\cline{2-7}
   	\multirow{2}{*}{~} &In-ADMM & 21&6.14/8 & 737.29 & $3.41\times 10^{-7}$ &$6.47\times10^{-7}$ \\
   	\hline
   		\multirow{2}{*}{$2^{-7}$} &Dual+SSN & 4 &13.25/15& 2710.77   & $3.41\times 10^{-7}$ & $6.47\times 10^{-7}$ \\
   		\cline{2-7}
   		\multirow{2}{*}{~} &In-ADMM & 20&5.85/8 & 6899.69 & $3.41\times 10^{-7}$ &$6.47\times10^{-7}$ \\
   		\hline
   		\multirow{2}{*}{$2^{-8}$} &Dual+SSN & 4 &14.25/17& 19112.53   & $3.41\times 10^{-7}$ & $6.47\times 10^{-7}$ \\
   		\cline{2-7}
   		\multirow{2}{*}{~} &In-ADMM & 17&6.11/8 &60794.04 & $3.41\times 10^{-7}$ &$6.47\times10^{-7}$ \\
   		\hline
   	\end{tabular}
   	\label{numerical_table3}
   \end{table}

    We notice that the outer iteration number of 'Dual+SSN' is much less than 'In-ADMM', which verified local superlinear convergence rate of SSN type method. Although each iteration of 'Dual+SSN' and 'In-ADMM' are both 2-nested layer, the main computational amount of 'Dual+SSN' is reflected in the process of solving linear equation w.r.t. preconditioner. That linear equation is solved by forward sweep and backward sweep which can also be treated as solving two discretized parabolic equations. That can help us interpret the total CPU time. 

\begin{figure}[H]
	\centering
	\begin{minipage}[t]{0.5\textwidth}
		\centering
		\includegraphics[width=6.5cm]{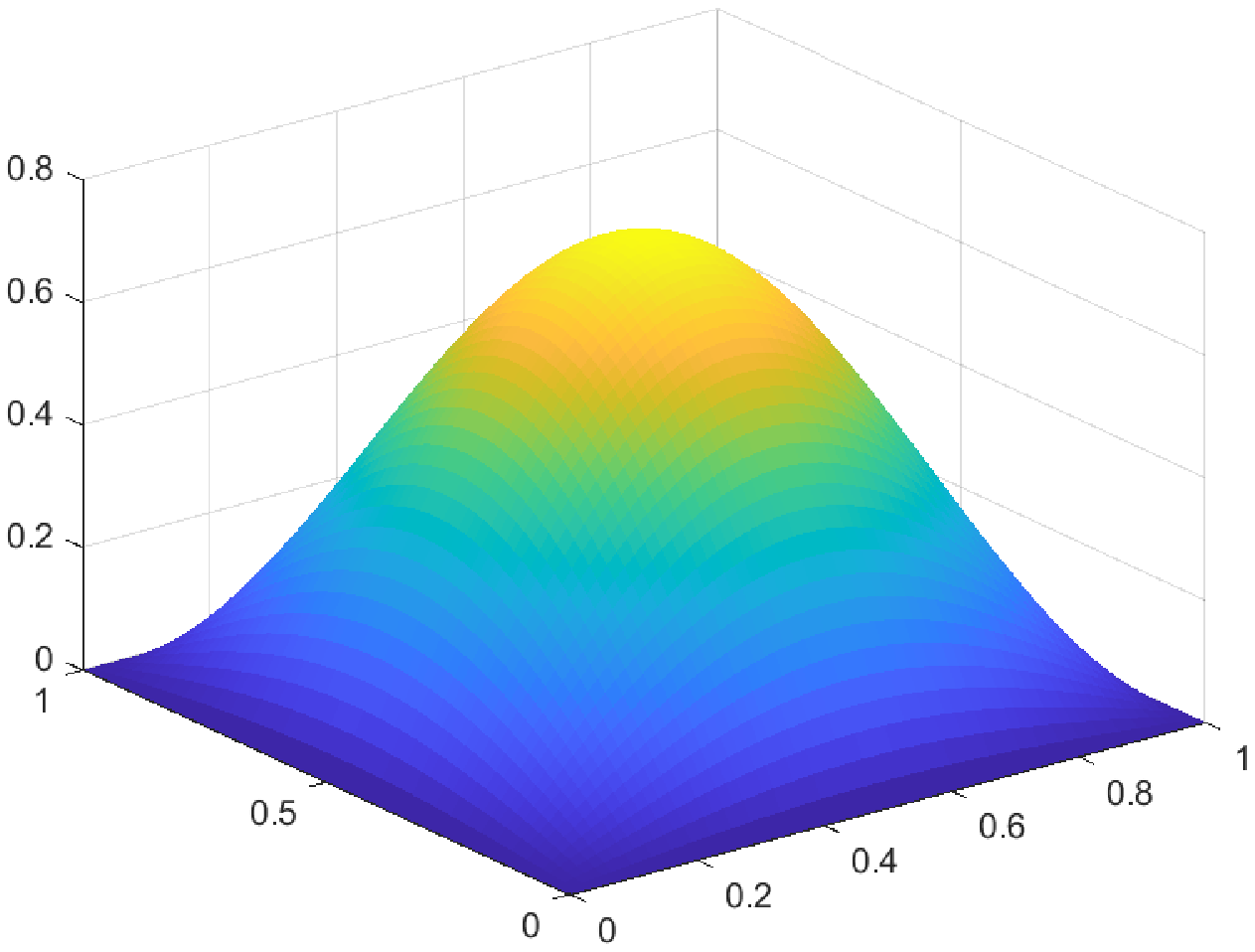}	
	\end{minipage}
	\begin{minipage}[t]{0.49\textwidth}
		\centering
		\includegraphics[width=6.5cm]{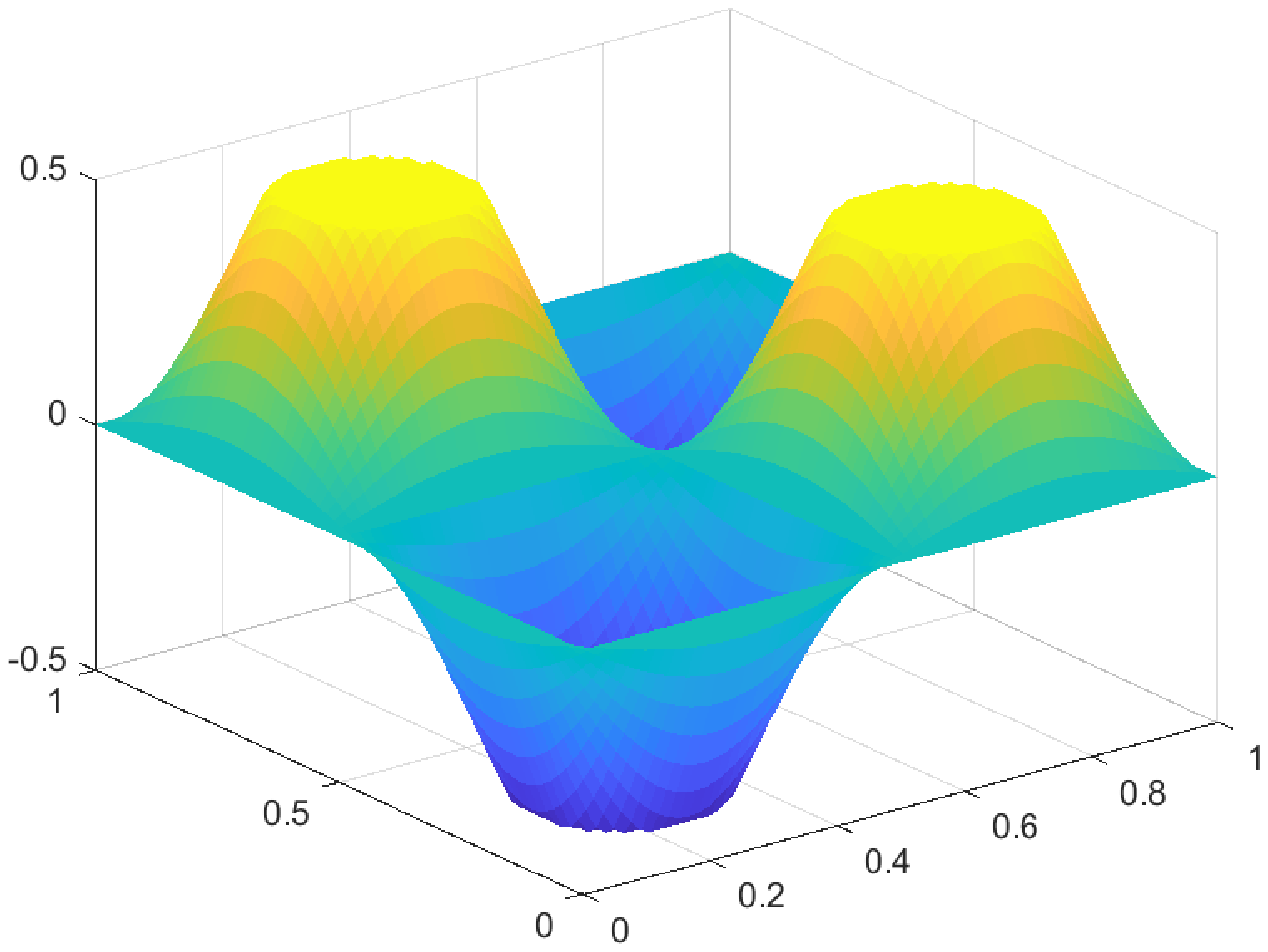}
	\end{minipage}
	\caption{Numercial solution (Dual+SSN) y (left) and u (right) at t=0.25 for Example 1}
	
\end{figure}

\begin{figure}[H]
	\centering
	\begin{minipage}[t]{0.5\textwidth}
		\centering
		\includegraphics[width=6.5cm]{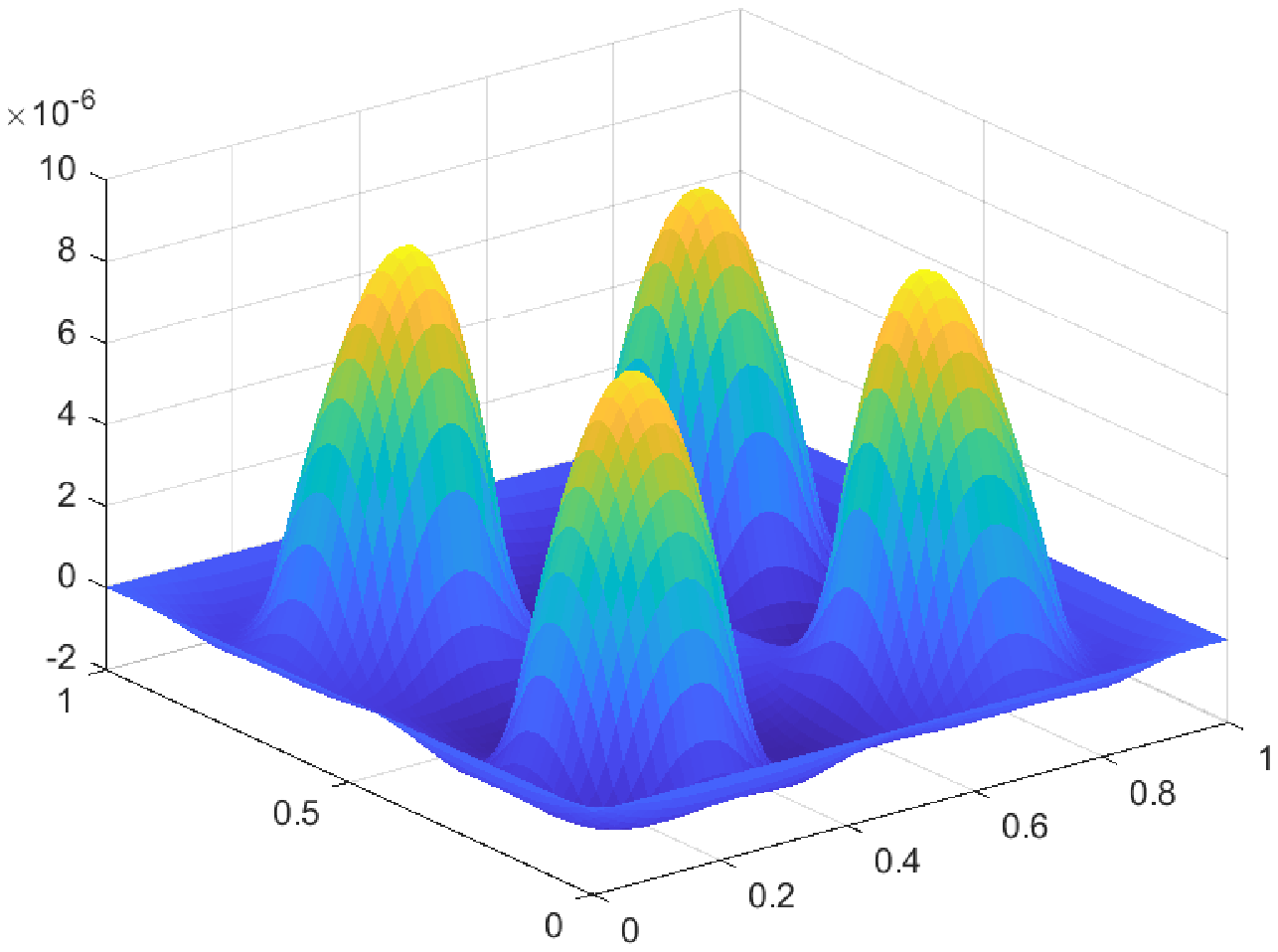}
	\end{minipage}
	\begin{minipage}[t]{0.49\textwidth}
		\centering
		\includegraphics[width=6.5cm]{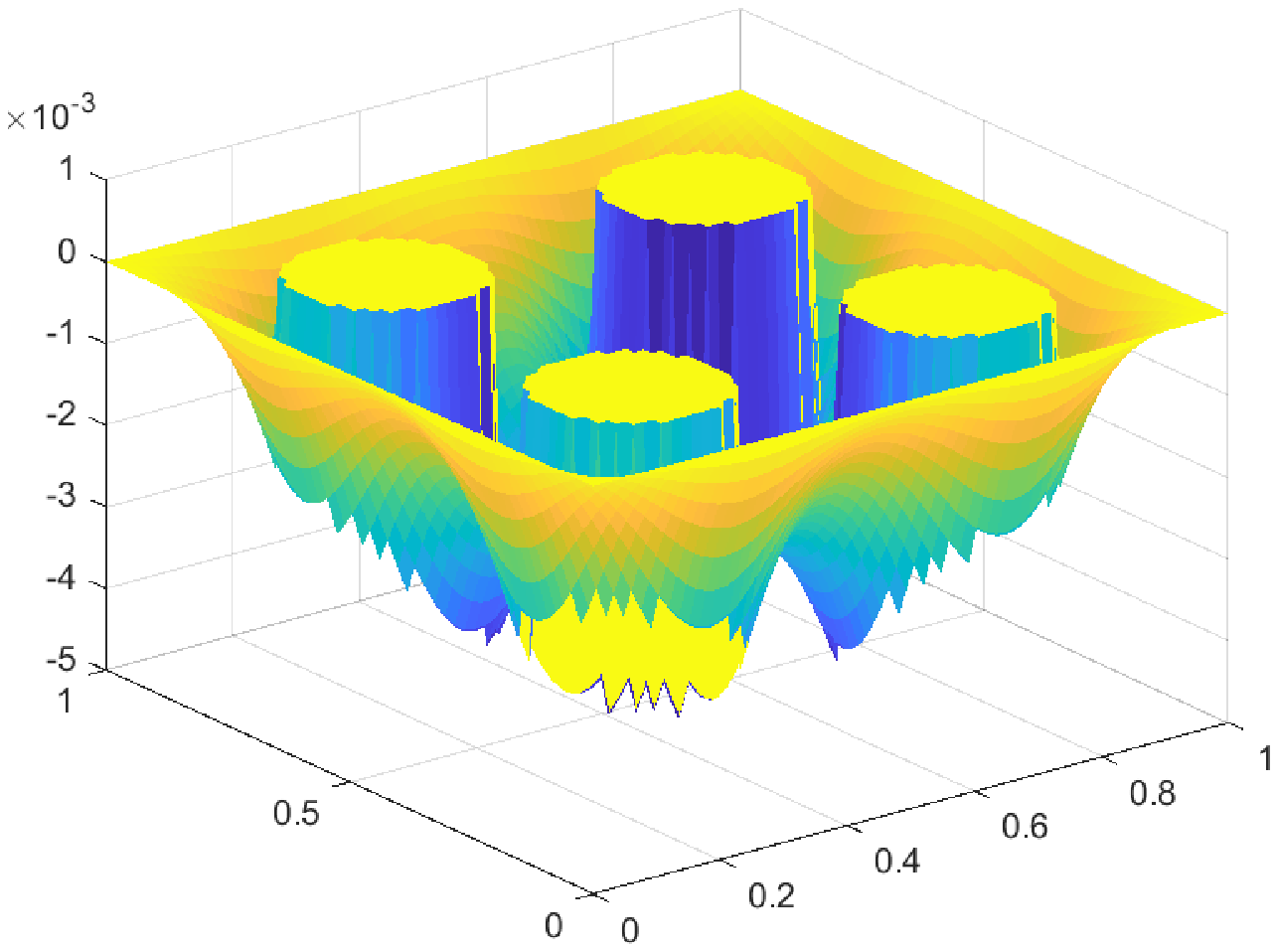}
	\end{minipage}
	\caption{Errors (Dual+SSN) $y-y^*$ (left) and $u-u^*$ (right) at t=0.25 for Example 1}
\end{figure}

   The following table verifies the accuracy of 'In-ADMM' and 'Dual+SSN', we observe that the error is dominated by the discretization error. Hence the accuracy of our designed algorithm 'Dual+SSN' can be guaranteed.
    \begin{table}[H]
   	\caption{Numerical errors comparsion of 'In-ADMM' and 'Dual+SSN' with $\gamma=10^{-5}$}
   	\centering
   	\begin{tabular}{| l |c | c| c|c | c|}
   		\hline 
   		error &Algorithm &$h=\Delta t=2^{-5}$&$h=\Delta t=2^{-6}$&$h=\Delta t=2^{-7}$&$h=\Delta t=2^{-8}$ \\
   		\hline
   		\multirow{2}{*}{$\|u-u^*\|_{L^2(\mathcal{O})}$} &Dual+SSN &$5.39\times10^{-3}$ &$1.37\times10^{-3}$ &$3.43\times10^{-4}$ & $8.57\times 10^{-5}$  \\
   		\cline{2-6}
   		\multirow{2}{*}{~} &In-ADMM &$5.39\times10^{-3}$ &$1.37\times10^{-3}$ &$3.57\times10^{-4}$ & $1.13\times 10^{-4}$  \\  		
   		\hline
   		\multirow{2}{*}{$\|y-y^*\|_{L^2(Q)}$} &Dual+SSN &$8.45\times10^{-6}$ &$2.15\times10^{-6}$ & $5.43\times 10^{-7}$ &$1.36\times10^{-7} $  \\
   		\cline{2-6}
   		\multirow{2}{*}{~} &In-ADMM &$8.47\times10^{-6}$ &$2.17\times10^{-6}$ & $5.81\times 10^{-7}$ &$1.91\times10^{-7} $\\
   		\hline
   	\end{tabular}
   	\label{numerical_error1}
   \end{table}
    
	\textbf{Example 2.} We consider another case where the control region $O$ is subset of domain $\Omega$, precisely $O=(0,0.25)\times(0,0.25)$ and $\Omega=(0,1)\times(0,1)$. Here we set $Q=\Omega\times(0,T)$ and $\mathcal{O}=O\times(0,T)$ and $T=1$. This problem is more general compared with Example 1 and its exact solution is unknown.
	\begin{equation*}
	\begin{aligned}
	\min\limits_{u \in \mathcal{C}, y \in L^{2}(Q)}\quad & \frac{1}{2} \iint_{Q}\left|y-y_{d}\right|^{2} d x d t+\frac{\gamma}{2} \iint_{\mathcal{O}}|u|^{2} d x d t \\
	\mbox{ s.t. }\quad & \left\{\begin{array}{ll}
	\frac{\partial y}{\partial t}-\Delta y+y=u\cdot\chi_{\mathcal{O}}, & \text { in } \Omega \times(0, T) \\
	y=0, & \text { on } \Gamma \times(0, T) \\
	y(0)=\sin(\pi x_1)\sin(\pi x_2)
	\end{array}\right.
	\end{aligned}
	\end{equation*}
	 The target function $y_d$ is specified by
	\begin{equation*}
	y_d=\exp(t)\sin(\pi x_1)\sin(\pi x_2)
	\end{equation*}
	and admissible set be
	\begin{equation*}
	\mathcal{C}=\{v|v\in L^2(\mathcal{O}), -300\leq v(t,x)\leq 300\text { a.e. in } Q \}
	\end{equation*}
	
	Firstly, we set regularization parameter $\gamma=10^{-3}$. We compare 'Dual+FRCG' with 'In-ADMM' in this case. The parameter $\beta$ appears in 'In-ADMM' is still set as $\beta=3$. From Table \ref{numerical_table2}, we observe that our method also behaves very efficient and robust for small control region case. Similar analysis as those we done for Example 1 can also be done for this example.

	\begin{table}[H]
		\caption{Numerical Comparsion of 'In-ADMM' and 'Dual+FRCG' when $\gamma=10^{-3}$}
		\centering
		\begin{tabular}{| l |c|c | c| c|c | c|}
			\hline 
			Mesh &Algorithm &Iter& Mean/Max CG &CPU Time(sec)&Obj &RelDis \\
			\hline
			\multirow{2}{*}{$2^{-4}$} &Dual+FRCG & 3 &\textemdash& 1.18   & $3.45\times 10^{-1}$ & $8.94\times 10^{-1}$ \\
			\cline{2-7}
			\multirow{2}{*}{~} &In-ADMM & 26&1/1 & 7.39  & $3.73\times 10^{-1}$ &$9.09\times10^{-1}$ \\
			\hline
			\multirow{2}{*}{$2^{-5}$} &Dual+FRCG & 3 &\textemdash& 3.84   & $3.64\times 10^{-1}$ & $9.2\times 10^{-1}$ \\
			\cline{2-7}
			\multirow{2}{*}{~} &In-ADMM & 26&1/1 & 22.09  & $3.78\times 10^{-1}$ &$9.27\times10^{-1}$ \\
			\hline
			\multirow{2}{*}{$2^{-6}$} &Dual+FRCG & 3 &\textemdash&25.97 &$3.73\times10^{-1}$ & $9.32\times10^{-1}$\\
			\cline{2-7}
			\multirow{2}{*}{~} &In-ADMM & 26&1/1 & 148.05  & $3.81\times 10^{-1}$ &$9.35\times10^{-1}$ \\
			\hline
			\multirow{2}{*}{$2^{-7}$} &Dual+FRCG & 3 &\textemdash&312.42 & $3.78\times 10^{-1}$ & $9.38\times10^{-1}$\\
			\cline{2-7}
			\multirow{2}{*}{~} &In-ADMM & 26&1/1 & 1809.19  & $3.82\times 10^{-1}$ &$9.40\times10^{-1}$ \\
			\hline
			\multirow{2}{*}{$2^{-8}$} &Dual+FRCG & 3 &\textemdash&2596.53   & $3.81\times 10^{-1}$ & $9.41\times10^{-1}$\\
			\cline{2-7}
			\multirow{2}{*}{~} &In-ADMM & 26&1/1 & 14821.61  & $3.96\times 10^{-1}$ &$9.75\times10^{-1}$ \\
			\hline
		\end{tabular}
	\label{numerical_table2}
	\end{table}
	
	\begin{figure}[H]
		\centering
		\begin{minipage}[t]{0.5\textwidth}
			\centering
			\includegraphics[width=6.5cm]{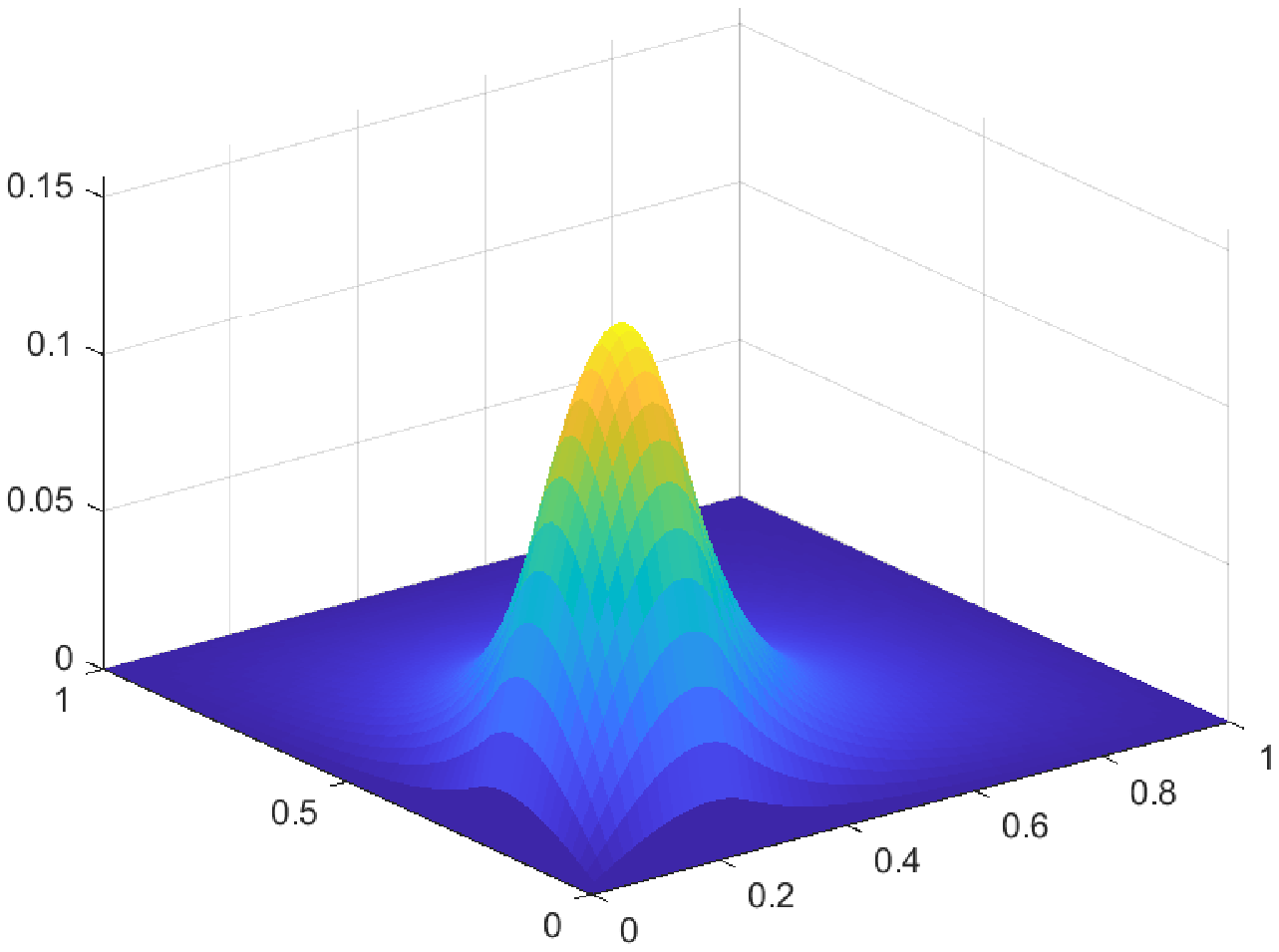}	
		\end{minipage}
		\begin{minipage}[t]{0.49\textwidth}
			\centering
			\includegraphics[width=6.5cm]{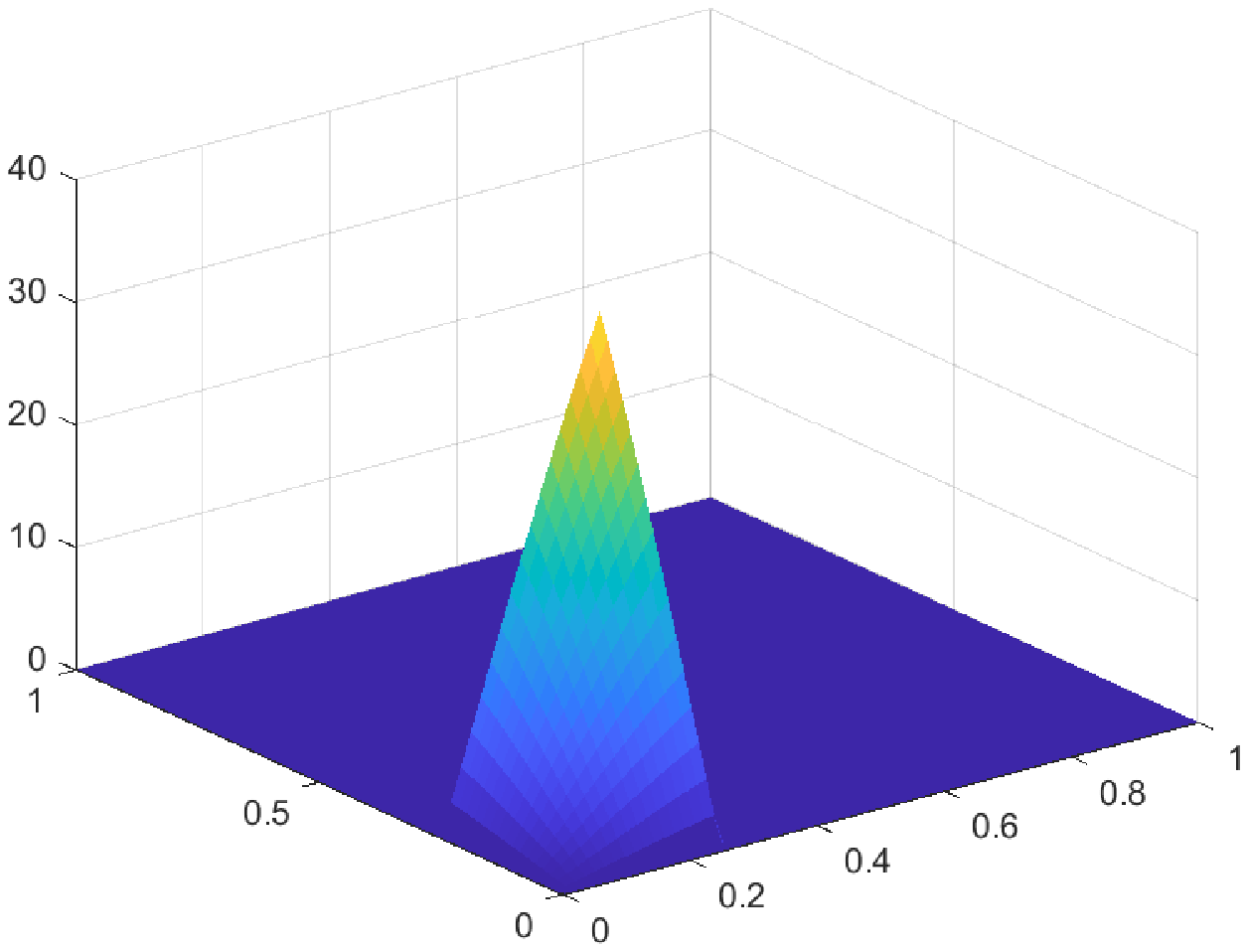}
		\end{minipage}
		\caption{Numercial solution (Dual+FRCG) y (left) and u (right) at t=0.5 for Example 2}
		
	\end{figure}
	 
 Then we set $\gamma=10^{-6}$ and compare 'Dual+SSN' with 'In-ADMM'. The parameter $\beta$ appears in 'In-ADMM' is still set as $\beta=3$. From Table \ref{numerical_table4}, we observe that 'In-ADMM' converges very slow this verifies its theoretical worst-case convergence rate while 'Dual+SSN' converges much faster. It seems that 'Dual+SSN' behaves much more efficient than 'In-ADMM' when the discretization is finer.
 \begin{table}[H]
 	\caption{Numerical Comparsion of 'In-ADMM' and 'Dual+SSN' when $\gamma=10^{-6}$}
 	\centering
 	\begin{tabular}{| l |c|c | c| c|c | c|}
 		\hline 
 		Mesh &Algorithm &Iter& Mean/Max CG &CPU Time(sec)&Obj &RelDis \\
 		\hline
 		\multirow{2}{*}{$2^{-4}$} &Dual+SSN & 6 &14.67/17& 13.64   & $2.66\times 10^{-1}$ & $6.62\times 10^{-1}$ \\
 		\cline{2-7}
 		\multirow{2}{*}{~} &In-ADMM & 90&3.05/5 & 60.38  & $2.66\times 10^{-1}$ & $6.62\times 10^{-1}$ \\
 		\hline
 		\multirow{2}{*}{$2^{-5}$} &Dual+SSN & 6 &17.67/21& 130.54   & $2.78\times 10^{-1}$ & $6.92\times 10^{-1}$ \\
 		\cline{2-7}
 		\multirow{2}{*}{~} &In-ADMM & 83&3.09/4 & 209.81  & $2.78\times 10^{-1}$ & $6.92\times 10^{-1}$ \\
 		\hline
 		\multirow{2}{*}{$2^{-6}$} &Dual+SSN & 7 &19.14/23& 1344.02   & $2.85\times 10^{-1}$ & $7.1\times 10^{-1}$ \\
 		\cline{2-7}
 		\multirow{2}{*}{~} &In-ADMM & 82&2.96/4 & 1384.56  & $2.85\times 10^{-1}$ & $7.1\times 10^{-1}$ \\
 		\hline
 		\multirow{2}{*}{$2^{-7}$} &Dual+SSN & 8 &19.75/25& 9014.77   & $2.89\times 10^{-1}$ & $7.2\times 10^{-1}$ \\
 		\cline{2-7}
 		\multirow{2}{*}{~} &In-ADMM & 83&2.94/3 & 16979.09  & $2.89\times 10^{-1}$ & $7.2\times 10^{-1}$ \\
 		\hline
 		\multirow{2}{*}{$2^{-8}$} &Dual+SSN & 7 &21.14/27& 55814.63   & $2.91\times 10^{-1}$ & $7.25\times 10^{-1}$ \\
 		\cline{2-7}
 		\multirow{2}{*}{~} &In-ADMM & 83&2.94/3 & 110283.29  & $2.91\times 10^{-1}$ & $7.25\times 10^{-1}$ \\
 		\hline
 	\end{tabular}
 	\label{numerical_table4}
 \end{table}
 
 \begin{figure}[H]
 	\centering
 	\begin{minipage}[t]{0.5\textwidth}
 		\centering
 		\includegraphics[width=6.5cm]{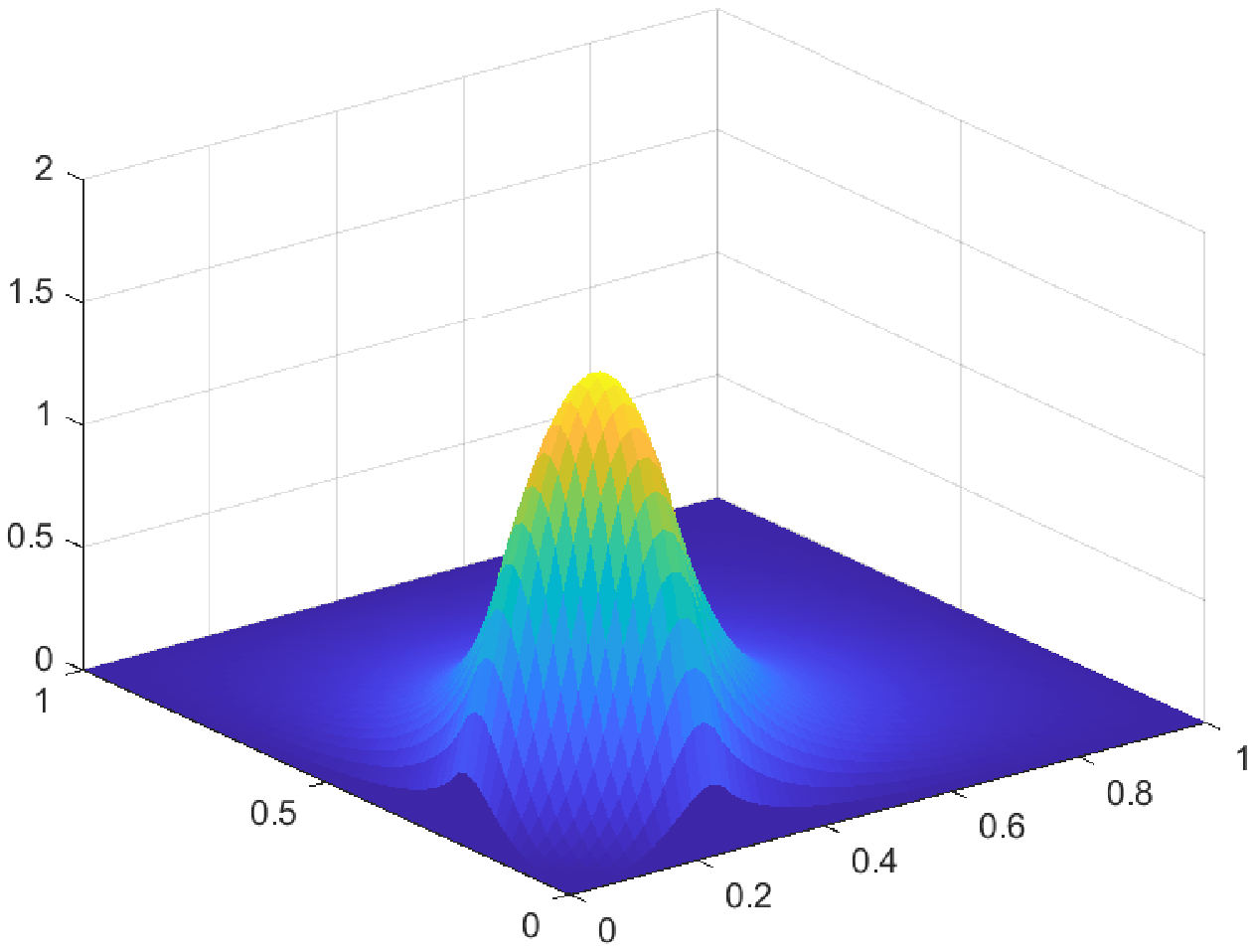}	
 	\end{minipage}
 	\begin{minipage}[t]{0.49\textwidth}
 		\centering
 		\includegraphics[width=6.5cm]{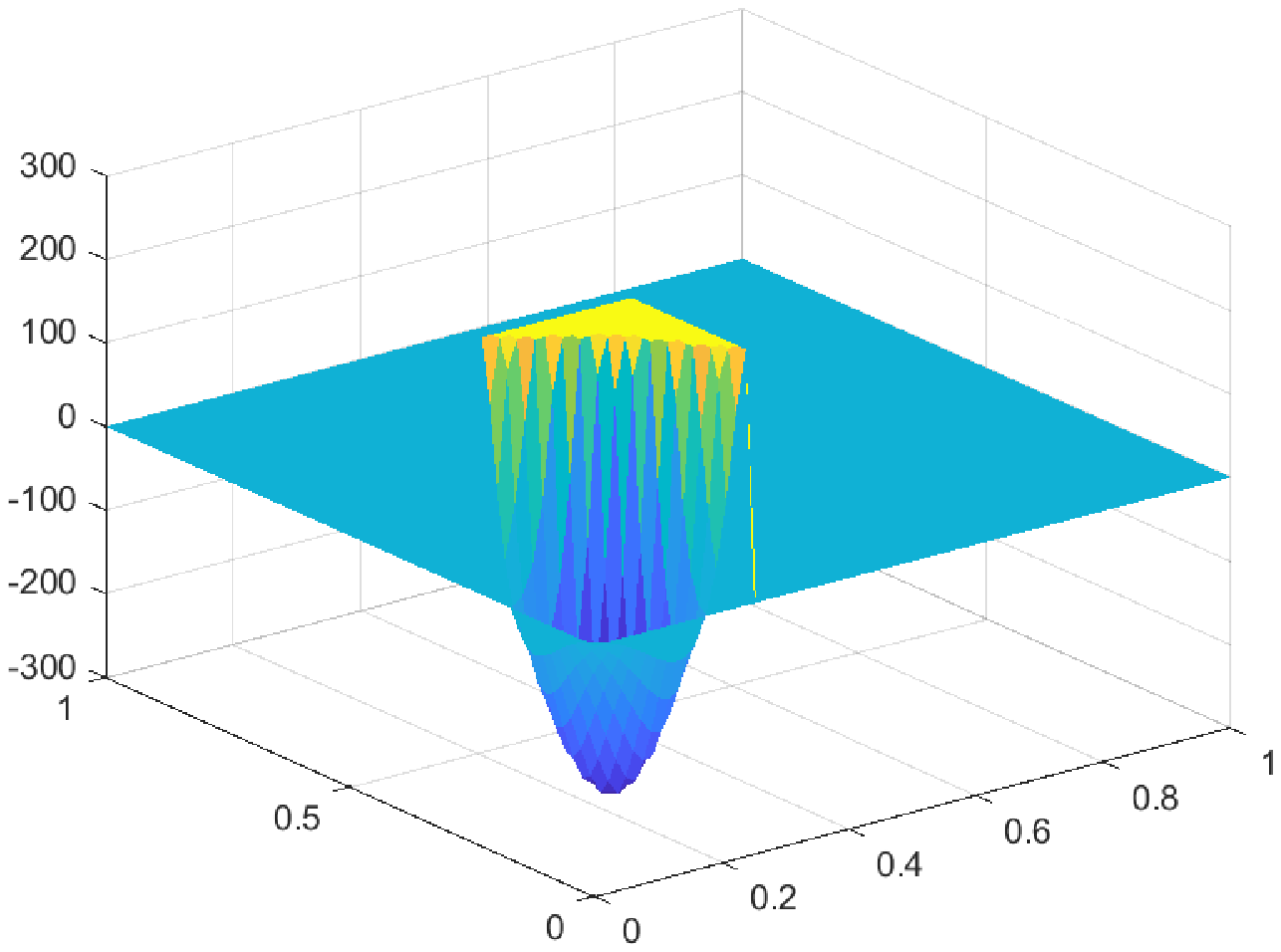}
 	\end{minipage}
 	\caption{Numercial solution (Dual+SSN) y (left) and u (right) at t=0.5 for Example 2}
 	
 \end{figure}

\textbf{Example 3.} Finally we test the second-order algorithm 'Dual+SSN' for elliptic optimal control problem with control constraints and compare it with the SSN method in \cite{MR3504560}. 

 We consider the following example given in \cite{MR2955288}:
\begin{equation*}
\begin{aligned}
	\min _{y \in H_{0}^{1}(\Omega), u \in \mathcal{C}} J(y, u)=& \frac{1}{2}\left\|y-y_{d}\right\|_{L^{2}(\Omega)}^{2}+\frac{\gamma}{2}\|u\|_{L^{2}(\Omega)}^{2} \\
	\text { s.t. } &\left\{\begin{array}{l}
		-\Delta y=u \text { in } \Omega, \\
		y=0 \text { on } \Gamma,
	\end{array}\right.
\end{aligned}
\end{equation*}

Let $\Omega=\left\{\left(x_{1}, x_{2}\right) \in \mathbb{R}^{2} \mid 0<x_{1}<1,0<x_{2}<1\right\}$ and the admissible set is specified as:
\begin{equation*}
\mathcal{C}=\left\{v \in L^{2}(\Omega) \mid -0.3 \leq u\left(x_{1}, x_{2}\right) \leq 1 \text { a.e. in } \Omega\right\}
\end{equation*}

The desired state is given by $y_{d}=4 \pi^{2} \gamma \sin \left(\pi x_{1}\right) \sin \left(\pi x_{2}\right)+y_{r}$. 

Here, the function $y_{r}$ denotes the solution to the following Possion equations:
\begin{equation*}
\begin{aligned}
-\Delta y_{r}&=r \quad\text { in } \Omega\\
 y_{r}&=0 \quad\text { on } \Gamma.
\end{aligned}
\end{equation*}
where $r=\min \left\{1, \max \left\{-0.3,2 \sin \left(\pi x_{1}\right) \sin \left(\pi x_{2}\right)\right\}\right\}$. It follows from the construction of $y_{d}$ and $r$ that $u^{*}:=r$ is the unique solution of this example.

To solve this problem, we can firstly derive its dual problem and then employ the SSN method for solving its discretized optimality condition obtained from dual problem. This procedure is totally similar to previous section \ref{SSN}, hence we omit specific detail.

For the numerical implementation of the SSN method, we follow the steps described in \cite{ MR3504560}. The initial values of the SSN method are set as $y=0,u=0,p=0$ and $\mu=0$, where $\mu=\mu_{a}+\mu_{b}$ with $\mu_{a}, \mu_{b}$ the Lagrange multipliers associated with the lower and upper bound of control constraints, as defined by equation (2.2) in \cite{ MR3504560}. 

We terminate SSN iterations when the nonlinear residual $F\left(u_{k} ; y_{k} ; p_{k} ; \mu_{k}\right) \leq 10^{-8}$( see (2.4) in \cite{ MR3504560}). We set $\gamma=10^{-4}$ in  and test various mesh sizes $h=2^{-i}$ with $i=4,5,6,7,8$. Besides in this example we terminate 'Dual+SSN' by (\ref{SSN_tol}) and cosntant is set as $tol=10^{-8}$. For both SSN type methods, we choose preconditioned GMRES method \cite{MR848568,MR1990645} for solving the obtained Newton equation at each step and tolerence is chosen as $10^{-8}$.

 Numerical results of the SSN in \cite{MR3504560} and 'Dual+SSN' iterative scheme are reported in the following table.

	\begin{table}[H]
	\caption{Numerical Comparsion of 'SSN' and 'Dual+SSN' }
	\centering
	\begin{tabular}{| l |c|c | c| c|c | c|}
		\hline 
		Mesh &Algorithm &Iter& Total GMRES &CPU Time(sec) &RelDis &$\|u-u^*\|_{L^2(\Omega)}$ \\
		\hline
		\multirow{2}{*}{$2^{-4}$} &Dual+SSN & 5&43&0.02 &$5.16\times 10^{-2}$  & $3.37\times10^{-4} $\\
		\cline{2-7}
		\multirow{2}{*}{~} &SSN & 5&41&0.08 &$5.16\times 10^{-2}$  & $3.44\times10^{-4}$\\
		\hline
		\multirow{2}{*}{$2^{-5}$} &Dual+SSN & 5&47&0.11 &$5.18\times 10^{-2}$  & $8.15\times10^{-5} $\\
	\cline{2-7}
	\multirow{2}{*}{~} &SSN & 5&40&0.13 &$5.18\times 10^{-2}$  & $8.15\times10^{-5}$\\
		\hline
		\multirow{2}{*}{$2^{-6}$} &Dual+SSN & 6&60&0.63 &$5.19\times 10^{-2}$  & $2.05\times10^{-5} $\\
		\cline{2-7}
		\multirow{2}{*}{~} &SSN & 6&47&0.61 &$5.19\times 10^{-2}$  & $2.2\times10^{-5}$\\
		\hline
		\multirow{2}{*}{$2^{-7}$} &Dual+SSN & 6&60&2.53 &$5.19\times 10^{-2}$  & $5.11\times10^{-6} $\\
		\cline{2-7}
		\multirow{2}{*}{~} &SSN & 6&43&2.21 &$5.19\times 10^{-2}$  & $8.03\times10^{-6}$\\
		\hline
		\multirow{2}{*}{$2^{-8}$} &Dual+SSN & 6&62&13.37 &$5.19\times 10^{-2}$  & $1.27\times10^{-6} $\\
		\cline{2-7}
		\multirow{2}{*}{~} &SSN & 5&35&9.77 &$5.19\times 10^{-2}$  & $2.37\times10^{-6}$\\
		\hline
	\end{tabular}
\label{numerical_table5}
\end{table}

\begin{figure}[H]
	\centering
	\begin{minipage}[t]{0.5\textwidth}
		\centering
		\includegraphics[width=6.5cm]{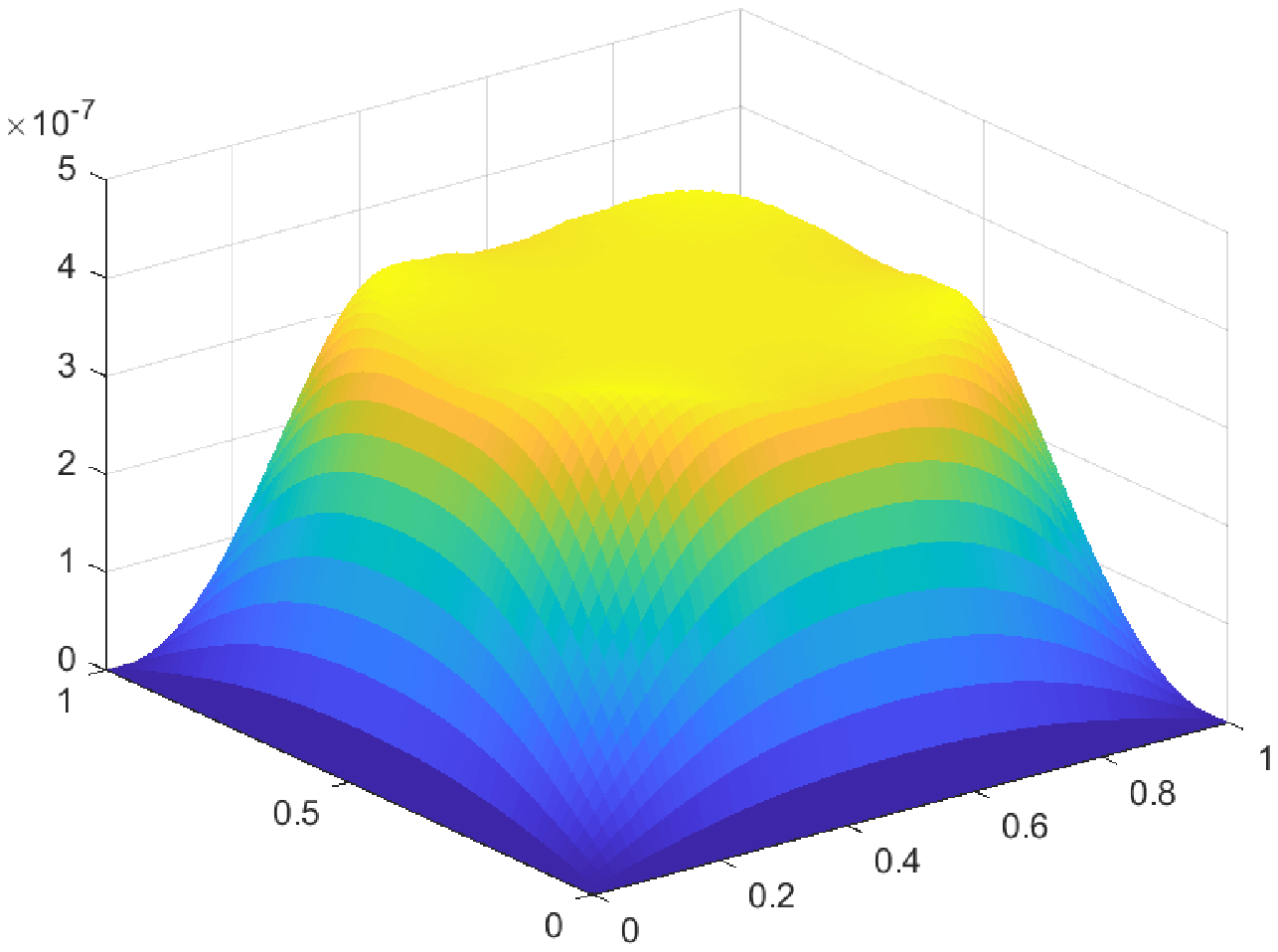}
	\end{minipage}
	\begin{minipage}[t]{0.49\textwidth}
		\centering
		\includegraphics[width=6.5cm]{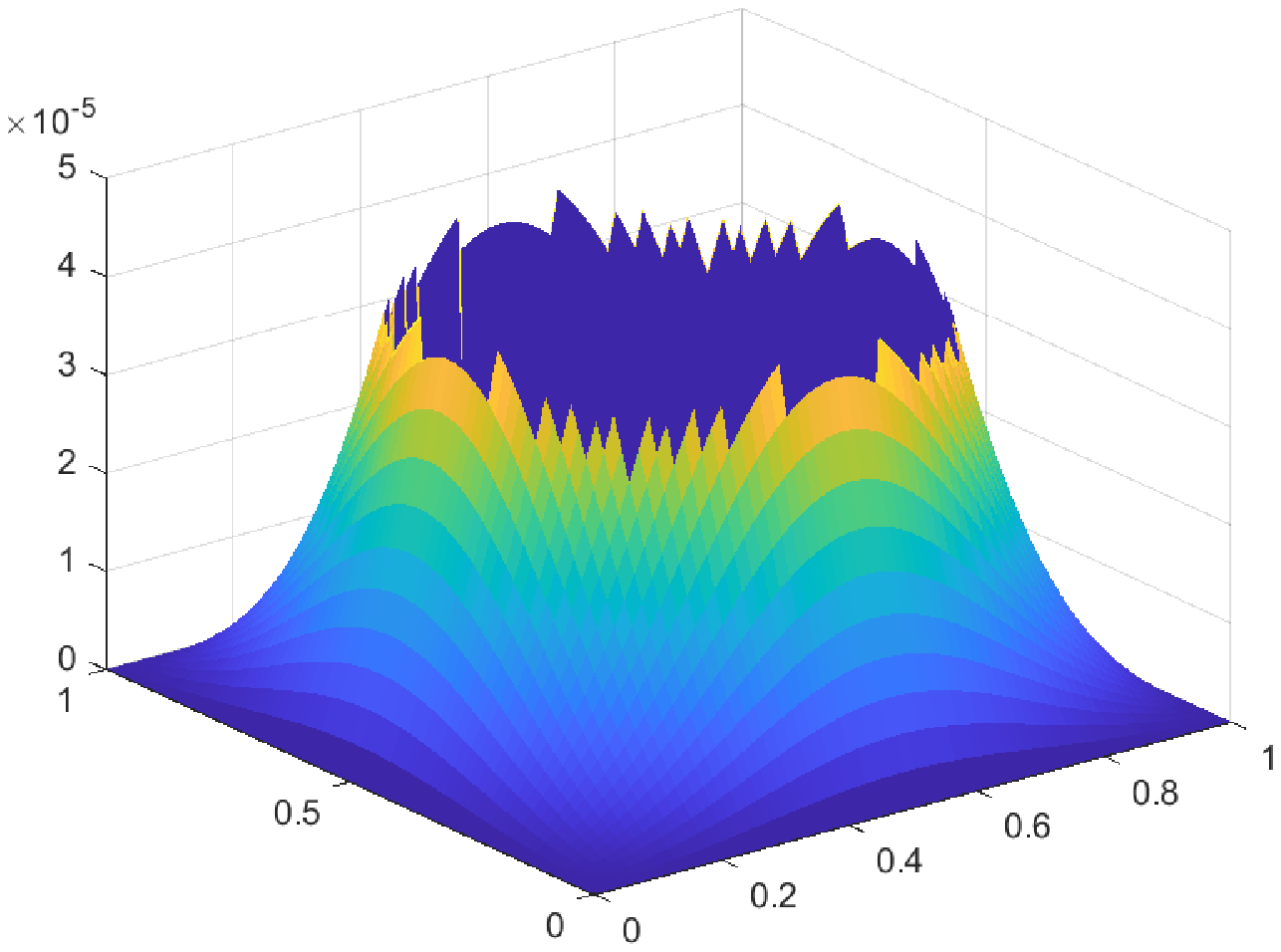}
	\end{minipage}
	\caption{Errors (Dual+SSN) $y-y^*$ (left) and $u-u^*$ (right) for Example 3}
\end{figure}

From Table \ref{numerical_table5}, we observe that the solution obtained by ‘Dual+SSN’ is even more accurate than SSN in \cite{ MR3504560}. For this example, our method requires more GMRES iterations to solve each Newton equation, but each Newton equation's scale is much smaller than that obtained by SSN. Besides, we notice that the outer iteration number is almost the same, thus total CPU time is of little difference. Hence, the 'Dual+SSN' is another an efficient method that can be used to solve elliptic optimal control problems.

\section{Conclusions}
In this paper, we focused on how to solve parabolic optimal control problems with control bounded constraints. Our algorithm design is totally based on dual problem that is different from ADMM type methods and SSN type methods. The dual problem derived by us is an unconstrained and first-order differentiable optimal control problem. And the control constraints occur in primal problem are embedded in the dual problem's objective functional. Indeed, the dual problem has better structure that helps us design more efficient numerical schemes. Besides we also specified the relationship between the solution of primal problem and the solution of dual problem. Our proposed method was first solve dual problem then obtain the solution of primal problem by the solution of dual problem. In order to solve dual problem numerically, we designed two easily implementable numerical schemes, precisely FRCG framework and SSN framework. For numerical discretization, we employed standard piecewise linear finite element method for space discretization and backward Euler finite difference method for time discretization. The resulting algorithms were validated to be numerically efficient by some preliminary numerical experiments.

Besides, this computational method advocated by us can be applied to other optimal control problems, including optimal control problems constrained by wave equations with control bounded constraints, elliptic optimal control problems with control constraints and linear diffusion-advection equations with control constraints etc. Conceptually our philosophy in algorithmic design can be extended to all optimal control problems with linear PDE constraints and control bounded constraints. In the future we will consider how to extend our method to solve more complicated control bounded optimal control problem e.g. dynamic system is described by nonlinear PDEs such as Navier-Stokes equations; sparsity optimal control problems formulated in non-reflexive Banach spaces, etc.    

\section*{Appendix.}\label{Appendix}
In this appendix, we provide the proof of Theorem \ref{th_spectral}

{\emph{Proof.}} Firstly we prove the lower bound: $\lambda\geq \frac{1}{2}$

We notice that there holds:
\begin{equation*}
\begin{aligned}
2C_k-\mathbb{C}_k&=2(\frac{\mathcal{M}_1\Pi_k}{\gamma}+\mathcal{K}\mathcal{M}^{-1}\mathcal{K}^{\top})-(\mathcal{K}+\frac{\mathcal{M}^{\frac{1}{2}}\mathcal{M}_1^{\frac{1}{2}}}{\sqrt{\gamma}}\Pi_k)\mathcal{M}^{-1}(\mathcal{K}+\frac{\mathcal{M}^{\frac{1}{2}}\mathcal{M}^{\frac{1}{2}}_1}{\sqrt{\gamma}}\Pi_k)^{\top}\\
&=\frac{\mathcal{M}_1\Pi_k}{\gamma}+\mathcal{K}\mathcal{M}^{-1}\mathcal{K}^{\top}-\frac{1}{\sqrt{\gamma}}(\mathcal{K}\mathcal{M}^{-\frac{1}{2}}\mathcal{M}_1^{\frac{1}{2}}\Pi_k+\mathcal{M}^{-\frac{1}{2}}\mathcal{M}_1^{\frac{1}{2}}\Pi_k\mathcal{K}^{\top})\\
&=(\mathcal{K}-\frac{\mathcal{M}^{\frac{1}{2}}\mathcal{M}_1^{\frac{1}{2}}}{\sqrt{\gamma}}\Pi_k)\mathcal{M}^{-1}(\mathcal{K}-\frac{\mathcal{M}^{\frac{1}{2}}\mathcal{M}^{\frac{1}{2}}_1}{\sqrt{\gamma}}\Pi_k)^{\top}
\end{aligned}
\end{equation*}

That means $2C_k-\mathbb{C}_k$ is positive semidefinite. Then there must holds:
\begin{equation*}
\frac{x^{\top}C_kx}{x^{\top}\mathbb{C}_kx}\geq \frac{1}{2} \quad \forall x
\end{equation*}

Meanwhile assume that $y$ is an eigenvector belongs to eigenvalue $\lambda$ w.r.t. matrix $\mathbb{C}_k^{-1}C_k$, then there holds:
\begin{equation*}
\lambda=\frac{y^{\top}C_ky}{y^{\top}\mathbb{C}_ky}
\end{equation*}

Hence there must hold $\lambda\geq \frac{1}{2}$

Next we should prove the upper bound. Still denotes $y$ as an eigenvector belong to eigenvalue $\lambda$ w.r.t. matrix $\mathbb{C}_k^{-1}C_k$

Let us denote $F$ as follows:
\begin{equation*}
F:=\frac{\mathcal{M}^{\frac{1}{2}}\mathcal{K}^{-1}\mathcal{M}_1^{\frac{1}{2}}\Pi_k}{\sqrt{\gamma}}
\end{equation*}

For eigenvalue $\lambda$ and eigenvector $y$, there holds $C_ky=\lambda\mathbb{C}_ky$. Then we have:
\begin{equation*}
\begin{aligned}
&(\frac{\mathcal{M}_1\Pi_k}{\gamma}+\mathcal{K}\mathcal{M}^{-1}\mathcal{K}^{\top})y=\lambda(\mathcal{K}+\frac{\mathcal{M}^{\frac{1}{2}}\mathcal{M}_1^{\frac{1}{2}}}{\sqrt{\gamma}}\Pi_k)\mathcal{M}^{-1}(\mathcal{K}+\frac{\mathcal{M}^{\frac{1}{2}}\mathcal{M}^{\frac{1}{2}}_1}{\sqrt{\gamma}}\Pi_k)^{\top}y\\
\iff
&(\frac{\mathcal{M}_1\Pi_k}{\gamma}+\mathcal{K}\mathcal{M}^{-1}\mathcal{K}^{\top})y=\lambda(\mathcal{K}\mathcal{M}^{-\frac{1}{2}}+\frac{\mathcal{M}_1^{\frac{1}{2}}}{\sqrt{\gamma}}\Pi_k)(\mathcal{K}\mathcal{M}^{-\frac{1}{2}}+\frac{\mathcal{M}_1^{\frac{1}{2}}}{\sqrt{\gamma}}\Pi_k)^{\top}y\\
\iff
&(\mathcal{K}\mathcal{M}^{-\frac{1}{2}})(I+FF^{\top})(\mathcal{K}\mathcal{M}^{-\frac{1}{2}})^{\top}y=\lambda (\mathcal{K}\mathcal{M}^{-\frac{1}{2}})(I+F)(I+F)^{\top}(\mathcal{K}\mathcal{M}^{-\frac{1}{2}})^{\top}y
\end{aligned}
\end{equation*}

Set $z=(I+F)^{\top}(\mathcal{K}\mathcal{M}^{-\frac{1}{2}})^{\top}y$, then we can derive the following equation:
\begin{equation*}
(I+F)^{-1}(I+FF^{\top})(I+F)^{-\top}z=\lambda z
\end{equation*}

Therefore we have:
\begin{equation*}
\begin{aligned}
\lambda&\leq \|(I+F)^{-1}(I+FF^{\top})(I+F)^{-\top}\|\\
&=\|(I+F)^{-1}\|^2+\|(I+F)^{-1}F\|^2\\
&=\|(I+F)^{-1}\|^2+\|I-(I+F)^{-1}\|^2\\
&\leq\|(I+F)^{-1}\|^2+(1+\|(I+F)^{-1}\|)^2\\
\end{aligned}
\end{equation*}

We set $\zeta=\|(I+F)^{-1}\|$, then the upper bounded is proved.

Finally we analyze the upper bound property when $\gamma\to 0^+$
This is equivalent to analyze $\|(I+F)^{-1}\|$. 

There holds:
\begin{equation*}
\begin{aligned}
\|(I+F)^{-1}\|&=\|\mathcal{M}^{\frac{1}{2}}(I+\frac{\mathcal{K}^{-1}\mathcal{M}_1^{\frac{1}{2}}\Pi_k\mathcal{M}^{\frac{1}{2}}}{\sqrt{\gamma}})^{-1}\mathcal{M}^{-\frac{1}{2}}\|\\
&\leq\operatorname{cond}(\mathcal{M}^{\frac{1}{2}})\|(I+\frac{\mathcal{K}^{-1}\mathcal{M}_1^{\frac{1}{2}}\Pi_k\mathcal{M}^{\frac{1}{2}}}{\sqrt{\gamma}})^{-1}\|
\end{aligned}
\end{equation*}

For the matrix $\mathcal{K}^{-1}\mathcal{M}_1^{\frac{1}{2}}\Pi_k\mathcal{M}^{\frac{1}{2}}$, we can do Jordan decomposition upon it. Thus we have the following equation:
\begin{equation*}
\mathcal{K}^{-1}\mathcal{M}_1^{\frac{1}{2}}\Pi_k\mathcal{M}^{\frac{1}{2}}=P\Lambda P^{-1}
\end{equation*}

Here $\Lambda$ is Jordan canonical form and suppose that 
\begin{equation*}
\Lambda=\operatorname{blkdiag}(\Lambda_1,\Lambda_2,\ldots,\Lambda_k,0,\ldots,0)
\end{equation*}
where each $\lambda_i (i=1,2,\ldots,k)$ denotes a nonzero Jordan block.

Since Lemma \ref{lemma_pd}, thus there must holds each eigenvalue of $\mathcal{K}^{-1}\mathcal{M}_1^{\frac{1}{2}}\Pi_k\mathcal{M}^{\frac{1}{2}}$ have strictly positive real part. Hence we can derive the following inequality:
\begin{equation*}
\begin{aligned}
\|(I+\frac{\mathcal{K}^{-1}\mathcal{M}_1^{\frac{1}{2}}\Pi_k\mathcal{M}^{\frac{1}{2}}}{\sqrt{\gamma}})^{-1}\|&=
\|(I+\frac{P\Lambda P^{-1}}{\sqrt{\gamma}})^{-1}\|\\
&=\|P(I+\frac{\Lambda}{\sqrt{\gamma}})^{-1}P^{-1}\|\\
&\leq \operatorname{cond}(P)\max\limits_{k}\{\|(I+\frac{\Lambda_k}{\sqrt{\gamma}})^{-1}\|,1\}
\end{aligned}
\end{equation*}

By the continuous property of norm, when $\gamma\to0^+$ we have $\|(I+\frac{\Lambda_k}{\sqrt{\gamma}})^{-1}\|\to0^+$. Thus we prove the conclusion. $\hfill\qedsymbol$

\clearpage
\bibliography{refee}

\end{document}